\documentclass{amsart}
\usepackage{amssymb}
\usepackage{amscd}

\theoremstyle{plain}
\newtheorem{thm}{Theorem}[section]
\newtheorem{lem}[thm]{Lemma}

\newtheorem{cor}[thm]{Corollary}

\theoremstyle{definition}
\newtheorem{dfn}[thm]{Definition}

\theoremstyle{remark}
\newtheorem{rem}[thm]{Remark}

\newcommand{\Z}{\mathbb{Z}}

\newcommand{\PP}{\mathbb{P}}
\newcommand{\OO}{\mathcal{O}}

\newcommand{\CC}{\mathcal{C}}
\newcommand{\XX}{\mathbb{X}}
\newcommand{\LL}{\mathbf{L}}
\newcommand{\pp}{\mathbf{p}}
\newcommand{\blambda}{\mbox{\boldmath $\lambda$}}

\DeclareMathOperator{\add}{add}
\DeclareMathOperator{\mmod}{mod}
\DeclareMathOperator{\Fac}{Fac}

\DeclareMathOperator{\lcm}{lcm}

\DeclareMathOperator{\Hom}{Hom}
\DeclareMathOperator{\ext}{ext}
\DeclareMathOperator{\Ker}{Ker}
\DeclareMathOperator{\Coker}{Coker}
\DeclareMathOperator{\Ima}{Im}

\DeclareMathOperator{\Aut}{Aut}
\DeclareMathOperator{\End}{End}
\DeclareMathOperator{\Ext}{Ext}

\DeclareMathOperator{\Rep}{Rep}

\DeclareMathOperator{\Floor}{Floor}

\DeclareMathOperator{\rank}{rank}

\begin{document}
\title{General sheaves over weighted projective lines}

\author{William Crawley-Boevey}
\address{Department of Pure Mathematics, University of Leeds, Leeds LS2 9JT, UK}
\email{w.crawley-boevey@leeds.ac.uk}

\subjclass[2000]{Primary 14H60; Secondary 16G20}



\keywords{Weighted projective line, coherent sheaf, general representation,
canonical decomposition}

\begin{abstract}
We develop a theory of general sheaves over weighted projective
lines. We define and study a canonical decomposition, analogous to Kac's
canonical decomposition for representations of quivers,
study subsheaves of a general sheaf, general ranks of morphisms,
and prove analogues of Schofield's results on general representations
of quivers. Using these, we give a recursive algorithm for
computing properties of general sheaves. Many of our results are
proved in a more abstract setting, involving a hereditary abelian category.
\end{abstract}
\maketitle

\section{Introduction}
In his work on representations of quivers, Kac~\cite{Kac1,Kac2,Kacm}
studied properties of the general representation of a given dimension
vector, for example he showed that any dimension vector $\alpha$
has a canonical decomposition, say $\alpha=\beta+\gamma+\dots$, such
that the general representation of dimension $\alpha$ is a direct sum of
representations of dimensions $\beta,\gamma,\dots$, each with trivial
endomorphism algebra. The theory of general representations was further
developed by Schofield \cite{S}, who found a recursive algorithm for
computing the canonical decomposition, as well as the dimensions of
Hom and Ext spaces between general representations, and possible
dimension vectors of subrepresentations of a general representation.
Later, Derksen and Weyman \cite{DWcan} and Schofield \cite{Sbq} found
efficient algorithms for computing the canonical decomposition.

The theory of general representations of quivers, and invariant theory
for representations of quivers, have been useful
in a number of areas, for example in Horn's problem on eigenvalues of
a sum of Hermitian matrices \cite{DWsemi}, in the theory of preprojective
algebras \cite{CBmm}, and in the study of vector bundles on
the projective plane \cite{Sprojplane}.
Therefore it would be nice to generalize this theory to quivers
with relations, or equivalently to the module varieties $\mmod_A^d$
classifying the $A$-module structures on a $d$-dimensional vector
space, for an arbitrary finitely generated algebra $A$.
Some progress is made in \cite{CBS}, but the theory is far from complete.

In this paper we instead study another hereditary situation, the
category of coherent sheaves on a weighted projective line in the
sense of Geigle and Lenzing \cite{GL}. This is an abelian category
which contains the category of vector bundles on $\PP^1$ equipped with
a quasi-parabolic structure of a certain type. We show that most of
the theory can be developed along similar lines in this setting.
We implicitly use the ideas of \cite{S} throughout.

Throughout $K$ is an algebraically closed field.
In the first part of the paper we work with a $K$-category $\CC$ which has
the following three properties.
\begin{itemize}
\item[(H)]$\CC$ is a hereditary abelian $K$-category with finite dimensional Hom and Ext spaces.
\item[(T)]For every object $X$ in $\CC$ there is a tilting object $T$ such that $X\in \Fac(T)$.
\item[(R)]For any objects $X$ and $Y$ the set $\Hom(X,Y)$ is partitioned into finitely many
locally closed subsets according to the rank of the morphism. The sets of monomorphisms and
epimorphisms are open.
\end{itemize}
Here $\CC$ being \emph{hereditary} means that $\Ext^n(X,Y)=0$ for $n\ge 2$.
Given an object $T$, one defines $\Fac(T)$ to be the set of objects which
are quotients of a finite direct sum of copies of $T$,
and $T$ is a \emph{tilting object} if $\Fac(T) = \{ X \mid \Ext^1(T,X)=0 \}$.
We consider $\Hom(X,Y)$ as affine space with the Zariski topology.
We call the class of an object $X$ in the Grothendieck group $K_0(\CC)$
its \emph{dimension type}, and denote it $[X]$.
The \emph{rank} of a morphism is the dimension type of its image.

Condition (H) ensures that the \emph{Euler form}
\begin{equation}\label{e:eulerform}
\langle X, Y \rangle = \dim\Hom(X,Y) - \dim\Ext^1(X,Y).
\end{equation}
depends only on the dimension types of $X$ and $Y$, so
it induces a bilinear form on $K_0(\CC)$.
We define
\[
K_0(\CC)_+ = \{ \alpha\in K_0(\CC) \mid \text{$\alpha=[X]$ for some object $X$} \}.
\]
Condition (T) implies that if $[X]=0$ then $X=0$,
since $X\in\Fac(T)$ for some tilting object $T$,
but then $\dim \Hom(T,X) = \langle T,X\rangle = \langle [T],0\rangle = 0$.
It follows that $K_0(\CC)_+$ is the positive cone for a partial ordering on $K_0(\CC)$.

In Section~\ref{s:variety} we introduce varieties $\Fac(T,\alpha)$ parametrizing
the objects in $\Fac(T)$ of dimension type $\alpha\in K_0(\CC)_+$, where $T$ is
a tilting object. These varieties depend on certain choices, and we study the
effect of changing the choices, or changing the tilting object in Section~\ref{s:change}.

We remark that, instead of using an infinite number of varieties to
parametrize the objects in $\CC$ of dimension type $\alpha$, we
should perhaps have used the language of stacks.

Another remark: our varieties $\Fac(T,\alpha)$ are essentially Quot schemes, except
that we avoid quotienting out by a certain group action. When $T$ is a tilting object,
however, these Quot schemes are particularly well behaved.

In Section~\ref{s:general} we define what it means for the
\emph{general object of dimension type $\alpha$} to have a given property $P$.
Following Kac, we show that any $\alpha\in K_0(\CC)_+$
has a \emph{canonical decomposition}
$\alpha=\beta_1+\beta_2+\dots+\beta_m$, unique up
to reordering, such that the general object of dimension
type $\alpha$ is a direct sum of indecomposable objects of
dimension types $\beta_1,\dots,\beta_m$.
We say that $\alpha\in K_0(\CC)_+$ is \emph{generally indecomposable}
if this canonical decomposition is trivial, that is, if the general object
of dimension type $\alpha$ is indecomposable, and we prove the following analogue
of Kac's result.

\begin{thm}\label{t:can}
The canonical decomposition
$\alpha=\beta_1+\beta_2+\dots+\beta_m$
is characterized by the fact that
the $\beta_i$ are generally indecomposable and
$\ext(\beta_i,\beta_j) = 0$ for $i\neq j$.
\end{thm}

Here $\hom(\beta,\gamma)$ and $\ext(\beta,\gamma)$ are the minimal dimensions
of $\Hom(X,Y)$ and $\Ext^1(X,Y)$ with $X$ of dimension $\beta$ and $Y$ of dimension
$\gamma$, or equivalently the general dimensions.

Given $\alpha,\gamma\in K_0(\CC)_+$ one defines $\gamma\hookrightarrow\alpha$ to
mean that the general object of dimension type $\alpha$ has a subobject
of dimension type $\gamma$. In Sections~\ref{s:subobjects} and \ref{s:generalrank}
we prove the analogues of Schofield's results on subobjects of a general object
and general ranks of morphisms.

\begin{thm}\label{t:subext}
Let $\alpha=\beta+\gamma$ with $\beta,\gamma\in K_0(\CC)_+$.
If $\ext(\gamma,\beta) = 0$, then $\gamma\hookrightarrow\alpha$.
If the field $K$ has characteristic 0, then the converse is also true.
\end{thm}

\begin{thm}\label{t:extform}
We have
\begin{align*}
\ext(\alpha,\beta)
&=
\min \{ - \langle \alpha-\eta,\beta-\eta\rangle
\mid \alpha-\eta \hookrightarrow \alpha, \eta\hookrightarrow \beta \}
\\
&= \min \{ - \langle \delta,\beta-\eta\rangle \mid
\delta \hookrightarrow \alpha, \eta\hookrightarrow \beta \}.
\end{align*}
\end{thm}

\begin{cor}\label{c:schur}
If the base field $K$ has characteristic 0, then the following are equivalent.

(1) $\alpha$ is generally indecomposable.

(2) The general object of dimension type $\alpha$ has one-dimensional endomorphism algebra.

(3) There is no nontrivial expression $\alpha=\beta+\gamma$ with $\beta\hookrightarrow\alpha$ and
$\gamma\hookrightarrow\alpha$.

(4) $\langle \beta,\alpha\rangle - \langle \alpha,\beta\rangle > 0$ for all $\beta\hookrightarrow\alpha$
with $\beta\neq 0,\alpha$.
\end{cor}

The equivalence of (1) and (2) was proved by Kac for representations
of quivers using completely different ideas. It is not clear how those
ideas can be adapted to this situation.

We now turn to weighted projective lines.
The category of coherent sheaves on a weighted projective line is well
known to have properties (H) and (T), and we show it has property (R)
in Section~\ref{s:wtd}. Therefore, the results above apply in this case.

Schofield's results easily give a recursive algorithm for computing
$\ext(\alpha,\beta)$ for representations of quivers.
Our analogous results Theorem~\ref{t:subext} and \ref{t:extform}
do not give an algorithm for weighted projective lines so easily,
but in Section~\ref{s:algorithm},
by exploiting the decomposition of a general sheaf into
the direct sum of a vector bundle and a torsion sheaf,
we are able to give an algorithm.

\begin{thm}\label{t:recalg}
For a weighted projective line, if the base field $K$ has characteristic zero,
then there is a recursive algorithm for computing $\ext(\alpha,\beta)$, and
hence also $\hom(\alpha,\beta)$, $\gamma\hookrightarrow\alpha$, and the
canonical decomposition.
\end{thm}

\begin{rem}
Recall that quasi-tilted algebras \cite{HRS} may be characterized as
the finite-dimensional algebras of global dimension at most 2, such
that any indecomposable module has projective dimension or injective
dimension at most 1. Equivalently, they are the endomorphism algebras
of tilting objects in $K$-categories $\CC$ with property (H).

Tilted algebras are the special case when $\CC$ is the category
of representations of a quiver. In \cite[\S 12.7]{CBS} it is shown how
the theory of general representations of quivers solves the problem of
classifying the irreducible components of module varieties for tilted
algebras $A$.

In the same way, our theory of general sheaves for
weighted projective lines should solve the problem of classifying the
irreducible components of module varieties for quasi-tilted algebras of
the form $A=\End(T)$, with $T$ a tilting sheaf for a weighted projective line.
\end{rem}

\begin{rem}
The theory of general representations of quivers is intimately related
to Schubert calculus. Exploiting the fact that Schubert calculus has been
developed in a characteristic-free way, this was used in \cite{CBschubert}
to extend Schofield's results to fields of positive characteristic.
It would be interesting to explore the relationship between the theory developed 
in this paper and quantum Schubert calculus \cite{B}. See also \cite{Belkale}.
\end{rem}

\section{A variety parametrizing objects}\label{s:variety}
In Sections \ref{s:variety}-\ref{s:generalrank} we work with a $K$-category $\CC$ with
properties (H), (T) and (R).
Let $T$ be a tilting object.  If $T''\in\add(T)$, that is,
$T''$ is isomorphic to a direct summand of a finite direct sum of copies of $T$,
then by a standard argument from tilting theory, the map $\Hom(T,T'')\to\Hom(T,X)$
induced by an epimorphism $f:T''\to X$ is onto if and only if $\Ker f\in\add(T)$.
In particular, for any $X\in \Fac(T)$, taking $f$ to be the universal morphism from
a direct sum of copies of $T$ to $X$, one sees that $X$ belongs to an exact sequence
\[
0\to T' \to T'' \to X \to 0
\]
with $T',T''\in\add(T)$.
If $Y$ is any object, then the universal extension
\[
0\to Y\to E\to T^{\oplus n}\to 0
\]
has the property that $\Ext^1(T,E)=0$, so that $E\in\Fac(T)$.
These exact sequences show that the dimension types of the non-isomorphic
indecomposable summands $T_i$ of $T$ generate $K_0(\CC)$, and they
are free generators since the argument of \cite[Lemma 4.1, Corollary 4.2]{HR}
shows that the matrix $\langle T_i, T_j \rangle$ is non-singular.
Thus also the Euler form is non-degenerate.

\begin{lem}
For any $\alpha\in K_0(\CC)$ there are objects
$T'_\alpha,T''_\alpha\in\add(T)$ with
\[
\alpha = [T''_\alpha] - [T'_\alpha]
\]
and such that any object $X$ in $\Fac(T)$ of dimension type $\alpha$
belongs to an exact sequence
\[
0\to T'_\alpha \to T''_\alpha \to X \to 0
\]
\end{lem}

\begin{proof}
If there is no such $X$, the existence of $T'_\alpha,T''_\alpha$
follows since the dimension types of the indecomposable summands
are a basis for $K_0(\CC)$.

Let $m = \langle T, \alpha\rangle$. We have $m\ge 0$, for if $X\in\Fac(T)$,
then $m = \dim\Hom(T,X)$. Take $T''_\alpha = T^{\oplus m}$. If $X\in\Fac(T)$,
then the universal map $T^{\oplus m}\to X$ is an epimorphism.
Its kernel is in $\add(T)$, and has dimension type $m[T] - \alpha$.
This determines the kernel up to isomorphism. We take it as $T'_\alpha$.
\end{proof}

Note that $T'_\alpha$ and $T''_\alpha$ are not uniquely determined in the lemma.
For example they can be replaced by $T'_\alpha\oplus S$ and $T''_\alpha\oplus S$
for any $S\in\add(T)$.

\begin{dfn}
Given $T$ and $\alpha\in K_0(\CC)$, we fix objects $T'_\alpha,T''_\alpha$ satisfying the conditions in the lemma.
We define $\Fac(T,\alpha)$ to be the set of monomorphisms $T''_\alpha\to T'_\alpha$.
It is an open subset of $\Hom(T''_\alpha,T'_\alpha)$, so a variety.
We define
\[
G(T,\alpha) = \Aut(T'_\alpha)\times\Aut(T''_\alpha).
\]
It is an algebraic group, and it acts on $\Fac(T,\alpha)$ via
\[
(f,g)\cdot\theta = g \theta f^{-1},
\quad
(f,g)\in \Aut(T'_\alpha)\times\Aut(T''_\alpha).
\]
\end{dfn}

Later we repeatedly use that $\Fac(T,\alpha)$ is an irreducible variety (if nonempty),
and for tangent space arguments we use that it is smooth.

If $\theta\in\Fac(T,\alpha)$ then there
is an exact sequence
\[
0\to T'_\alpha \xrightarrow{\theta} T''_\alpha \xrightarrow{p_\theta} \Coker\theta\to 0.
\]
so $\Coker\theta$ has dimension type $\alpha$. By hypothesis any object in $\Fac(T)$
of dimension type $\alpha$ is isomorphic to $\Coker\theta$ for some $\theta\in \Fac(T,\alpha)$.

For ease of notation we write $\Fac(T,\alpha)\times\Fac(T,\beta)$ as $\Fac(T,\alpha,\beta)$.
Let $C(T,\alpha,\beta)$ be the variety of tuples
\[
((\theta,\phi),f,g)\in \Fac(T,\alpha,\beta)\times\Hom(T'_\alpha,T'_\beta)\times\Hom(T''_\alpha,T''_\beta)
\]
with $g\theta = \phi f$, so forming a commutative square
\[
\begin{CD}
T'_\alpha @>\theta>> T''_\alpha \\
@Vf VV @Vg VV \\
T'_\beta @>\phi>> T''_\beta.
\end{CD}
\]
For any such commutative square there is a unique $h$ giving a morphism
of exact sequences
\begin{equation}\label{e:sesmorph}
\begin{CD}
0 @>>> T'_\alpha @>\theta>> T''_\alpha @>p_\theta >> \Coker \theta @>>> 0 \\
& & @Vf VV @Vg VV @Vh VV \\
0 @>>> T'_\beta @>\phi>> T''_\beta @>p_\phi >> \Coker \phi @>>> 0.
\end{CD}
\end{equation}

Conversely, given $\theta$ and $\phi$, any morphism $h$ arises in this way
for some $f,g$. More precisely, letting
\[
\pi:C(T,\alpha,\beta)\to \Fac(T,\alpha,\beta)
\]
be the projection, the following is straightforward.

\begin{lem}\label{l:homlifts}
If $(\theta,\phi)\in \Fac(T,\alpha,\beta)$ then there
is an exact sequence of vector spaces
\begin{equation}\label{e:pifibses}
0\to \Hom(T''_\alpha,T'_\beta) \to \pi^{-1}(\theta,\phi) \xrightarrow{\kappa} \Hom(\Coker\theta,\Coker\phi) \to 0,
\end{equation}
where $\kappa((\theta,\phi),f,g)$ is the unique $h$ in \textup{(\ref{e:sesmorph})}.
\end{lem}

\begin{lem}
The function $(\theta,\phi)\mapsto\dim\Hom(\Coker\theta,\Coker\phi)$
is upper semicontinuous on $\Fac(T,\alpha,\beta)$
\end{lem}

\begin{proof}
Let $s: \Fac(T,\alpha,\beta) \to C(T,\alpha,\beta)$ be the zero section
sending $(\theta,\phi)$ to $((\theta,\phi),0,0)$. Since the fibres of $\pi$,
are vector spaces, they intersect the image of the zero section.
The result thus follows from the Upper Semicontinuity Theorem applied to $\pi$.
\end{proof}

Similarly, for any fixed object $X$, the functions $\theta\mapsto\dim\Hom(X,\Coker\theta)$
and $\theta\mapsto\dim\Hom(\Coker\theta,X)$ are upper semicontinuous on $\Fac(T,\theta)$.
Also, thanks to the Euler form, the analogous functions with $\dim\Hom$
replaced by $\dim\Ext^1$ are upper semicontinuous.

\begin{lem}
There is a 1:1 correspondence between the orbits of $G(T,\alpha)$ on $\Fac(T,\alpha)$
and the isomorphism classes of objects in $\Fac(T)$ of dimension type~$\alpha$.
\end{lem}

\begin{proof}
We need to show that if $\theta,\phi\in\Fac(T,\alpha)$ and there
is an isomorphism $h:\Coker\theta\to\Coker\phi$, then $\theta$
and $\phi$ are in the same orbit under $G(T,\alpha)$.

One can decompose $T''_\alpha = X_0 \oplus X_1$ such that
$p_\theta|_{X_0}$ is right minimal and $p_\theta|_{X_1}=0$,
and $T''_\alpha = Y_0 \oplus Y_1$ such that $p_\phi|_{X_0}$
is right minimal and $p_\phi|_{X_1}=0$, see for example \cite{AS}.
By Lemma~\ref{l:homlifts}, $h$ lifts to an endomorphism of
$T''_\alpha$, and if $a$ is the component $X_0\to Y_0$
then $p_\phi|_{Y_0}a = h p_\theta|_{X_0}$.
Similarly $h^{-1}$ lifts to an endomorphism of
$T''_\alpha$, and if $b$ is the component $Y_0\to X_0$,
then $p_\theta|_{X_0} b = h^{-1} p_\phi|_{Y_0}$.
Thus by minimality $ab$ and $ba$ are automorphisms,
so $a$ and $b$ are isomorphisms. Now by the Krull-Remak-Schmidt Theorem,
$X_1\cong Y_1$, and taking the morphism $X_0\oplus X_1\to Y_0\oplus Y_1$
with components $a$ and an arbitrary isomorphism $X_1\to Y_1$, we
obtain an automorphism $g$ of $T''_\alpha$ with $p_\phi g = hp_\theta$.
There is now a unique morphism $f$, necessarily an automorphism,
completing the diagram (\ref{e:sesmorph}).
The group element $(f,g)$ acting on $\theta\in\Fac(T,\alpha)$ then gives $\phi$.
\end{proof}

Later we will need the following.

\begin{lem}
If $((\theta,\phi),f,g)\in C(T,\alpha,\beta)$, and $h$ is as in \textup{(\ref{e:sesmorph})},
then there are exact sequences
\begin{equation}\label{e:gphiseq}
0 \to \Ker(hp_\theta) \to T''_\alpha\oplus T'_\beta
\xrightarrow{(\begin{smallmatrix} g & \phi \end{smallmatrix})}
T''_\beta \to \Coker h \to 0
\end{equation}
and
\begin{equation}\label{e:kerhses}
0\to T'_\alpha\to \Ker(hp_\theta) \to \Ker h \to 0.
\end{equation}
\end{lem}

\begin{proof}
Straightforward.
\end{proof}

\begin{lem}
Given a tilting object $T$ and $\alpha,\beta\in K_0(\CC)$, there
are only finitely many possible ranks of morphisms from
an object of dimension $\alpha$ in $\Fac(T)$ to one
of dimension $\beta$ in $\Fac(T)$.
\end{lem}

\begin{proof}
By (\ref{e:gphiseq}), the rank of $h$ is determined by that
of $(g\ \phi)$, and by condition (R) there are only
finitely many possibilities for this.
\end{proof}

\section{Changing the tilting object or the presentation}\label{s:change}
Suppose $S$ and $T$ are tilting objects with $\Fac(S)\subseteq \Fac(T)$.
Given $\alpha$, one can consider $\Fac(S,\alpha)$ and $\Fac(T,\alpha)$,
say given by objects $S',S''$ and $T',T''$. (We drop the subscript $\alpha$
in this section, to simplify notation.)
In this section we relate the variety $\Fac(S,\alpha)$ and
\[
\Fac(T,\alpha)_S = \{ \theta\in \Fac(T,\alpha) \mid \Coker\theta\in\Fac(S)\}.
\]
As a special case this includes the possibility that $S=T$, but we choose
different objects of $\add(T)$ in the definition of $\Fac(T,\alpha)$.

Since $\Coker\theta\in\Fac(S)$ if and only if $\Ext^1(S,\Coker\theta)=0$,
upper semicontinuity of $\dim\Ext^1(S,\Coker\theta)$ implies that
$\Fac(T,\alpha)_S$ is an open subset of $\Fac(T,\alpha)$.

Let $E$ be the variety of commutative squares
\[
\begin{CD}
T' @>\theta>> T''\\\
@VfVV @VgVV \\
S' @>\phi>> S''
\end{CD}
\]
such that $\theta$ and $\phi$ are monomorphisms and
such that the sequence
\[
0\to T'\to S'\oplus T'' \to S''\to 0
\]
is exact. This last condition is open since
by dimensions one only needs to ensure that
$T'\to S'\oplus T''$ is a monomorphism and
$S'\oplus T''\to S''$ is an epimorphism.

We consider elements of $E$ as endomorphisms of
$T'\oplus T''\oplus S'\oplus S''$ of shape
\[
\begin{pmatrix}
0 & 0 & 0 & 0 \\
\theta & 0 & 0 & 0 \\
f & 0 & 0 & 0 \\
0 & -g & \phi & 0
\end{pmatrix}
\]
and square zero.
We consider the action by conjugation of the group $G$ of
automorphisms of $T'\oplus T''\oplus S'\oplus S''$ of the form
\[
\begin{pmatrix}
a & 0 & 0 & 0 \\
0 & b & 0 & 0 \\
0 & e & c & 0 \\
0 & 0 & 0 & d
\end{pmatrix}.
\]
The inverse matrix is
\[
\begin{pmatrix}
a^{-1} & 0 & 0 & 0 \\
0 & b^{-1} & 0 & 0 \\
0 & -c^{-1}eb^{-1} & c^{-1} & 0 \\
0 & 0 & 0 & d^{-1}
\end{pmatrix}.
\]
Thus the action is given by
\[
(a,b,c,d,e)\cdot (\theta,\phi,f,g) =
(b\theta a^{-1},d\phi c^{-1}, cfa^{-1}+e\theta a^{-1},dgb^{-1}+d\phi c^{-1}eb^{-1})
\]
The group $G$ has natural homomorphisms to $G(S,\alpha)$
and $G(T,\alpha)$, so it acts on $\Fac(S,\alpha)$ and $\Fac(T,\alpha)$.

\begin{lem}
The projection $E\to \Fac(S,\alpha)$ is smooth and onto, and it is
equivariant for the action of $G$.
The fibres are orbits for the kernel of the group homomorphism $G\to G(S,\alpha)$.
The inverse image of a $G(S,\alpha)$-orbit in $\Fac(S,\alpha)$ is a $G$-orbit in $E$.
\end{lem}

\begin{proof}
Inside the product $\Hom(T'',S'')\times \Fac(S,\alpha)$ look at the
open subset $P$ given by those pairs $(g,\phi)$
such that the map $\Hom(T,T''\oplus S')\to \Hom(T,S'')$ is onto.
Because of the exact sequence
\[
0\to \Hom(T,S')\to \Hom(T,S'')\to \Hom(T,\Coker\phi)\to 0
\]
it is equivalent that the map $\Hom(T,T'')\to \Hom(T,\Coker\phi)$ is onto.
Thus it is equivalent that the kernel of the map $T''\to\Coker\phi$ is
in $\add(T)$ (so isomorphic to $T'$).

There is a projection map $E\to P$, for the exact
sequence defining $E$ clearly shows that
$\Hom(T,T''\oplus S')\to \Hom(T,S'')$ is onto.

The fibre over $(g,\phi)$ is determined by an isomorphism from $T'$ to the
kernel of the map $T''\to\Coker\phi$.

Thus $E\to\Fac(S,\alpha)$ is smooth and onto.
Equivariant for action of $G(S,\alpha)$.
Constant on orbits of $G(T,\alpha)$.

Now suppose that $(\theta,\phi,f,g)$ and $(\theta',\phi,f',g')$ are in $E$ (same $\phi$).
There are induced isomorphisms $\Coker\theta\to\Coker\phi$ and $\Coker\theta'\to\Coker\phi$.
Let $\Coker\theta'\to\Coker\theta$ be the corresponding isomorphism.
Any such isomorphism lifts to an element $(a,b)\in \Aut(T')\times\Aut(T'') = G(T,\alpha)$.
Thus $\theta a = b\theta'$, and the induced map on cokernels is as given.
Thus the morphisms $\Coker\theta'\to\Coker\phi$ induced by $(f',g')$ and $(fa,gb)$ are
the same. Thus $(f'-fa,g'-gb)$ induces the zero map on cokernels.
Thus $g'-gb = \phi e$ for some (uniquely determined) morphism $e\in \Hom(T'',S')$.
Then
$\phi (f'-fa) = \phi f' - \phi fa = g'\theta' - g\theta a  = (gb+\phi e)\theta' - g\theta a = \phi e\theta'$
and so $f'-fa = e\theta' = eb^{-1}\theta a$.
Thus
\[
(\theta',\phi,f',g') = (a^{-1},b^{-1},1,1,eb^{-1})\cdot(\theta,\phi,f,g).
\]
\end{proof}

\begin{lem}
The projection $E\to \Fac(T,\alpha)$ is smooth and has image $\Fac(T,\alpha)_S$.
It is $G$-equivariant. The fibres are orbits for the kernel of
the group homomorphism $G\to G(T,\alpha)$. The inverse image of a $G(T,\alpha)$-orbit
in $\Fac(T,\alpha)$ is a $G$-orbit in $E$.
\end{lem}

\begin{proof}
Image is contained in $\Fac(T,\alpha)_S$
since $\Coker\theta\cong\Coker\phi\in\Fac(S)$.

Image is equal to $\Fac(T,\alpha)_S$ since if
$\theta\in\Fac(T,\alpha)_S$ then there is some $\phi\in\Fac(S,\alpha)$
with $\Coker\phi\cong\Coker\theta$.
Then an isomorphism
$\Coker\theta\cong\Coker\phi$ lifts to a morphism $g:T''\to S''$ since $\Ext^1(T'',S')=0$.
Then get induced morphism $f:T'\to S'$.

Fibre is orbit for the kernel of
the group homomorphism $G\to G(S,\alpha)$.
Suppose $(\theta,\phi,f,g)$ and $(\theta,\phi',f',g')$ are
both in $E$.
There are induced isomorphisms $\Coker\theta\to\Coker\phi$
and $\Coker\theta\to\Coker\phi'$.
Consider the corresponding isomorphism $\Coker\phi\to\Coker\phi'$.
Lift to $(c,d)\in\Aut(S')\times \Aut(S'')=G(S,\alpha)$, so $d\phi = \phi'c$.
Then $(f',g')$ and $(cf,dg)$ induce the same morphism $\Coker\theta\to\Coker\phi'$.
Thus $g'-dg = \phi'e$ for some $e\in\Hom(T'',S')$. etc.

Consider the open set $Q$ of pairs $(f,\theta)\in \Hom(T',S')\times\Fac(T,\alpha)_S$
such that the map
\[
\Ext^1(S,T')\to \Ext^1(S,T''\oplus S')
\]
induced by $T'\to T''\oplus S'$ is onto. It is equivalent that the
cokernel of $T'\to T''\oplus S'$ is in $\add(S)$, so equivalent
that the cokernel is isomorphic to $S''$.

Consider the trivial bundle $\Hom(T'',S'')\times \Hom(S',S'')\times Q\to Q$.
The subset $R$ consisting of elements $(g,\phi,(f,\theta))$ such that
$\phi f = g\theta$ is a vector sub-bundle. Thus the projection $R\to Q$
is smooth. Inside $R$, consider the subset for which the sequence
\[
0\to T'\to S'\oplus T'' \to S''\to 0
\]
(which we know has composite zero) is exact.
This defines an open subset of $R$.
This open set is identified with $E$. Thus $E\to Q$ is smooth.
Hence $E\to \Fac(T,\alpha)_S$ is smooth.
\end{proof}

Combining the above we get the following.

\begin{thm}\label{t:change}
There is a 1-1 correspondence between
$G(T,\alpha)$-stable open subsets of $\Fac(T,\alpha)_S$
and $G(S,\alpha)$-stable open subsets of $\Fac(S,\alpha)$,
corresponding open subsets being determined by the possible
isomorphism classes of objects they classify.
\end{thm}

\section{Properties of the general object}\label{s:general}
Let $P$ be a property of objects in $\CC$ of dimension type $\alpha$, depending only
on the isomorphism class of the object.

\begin{dfn}
Given $\alpha\in K_0(\CC)_+$, we say that
the \emph{general object of dimension type $\alpha$ has property $P$}
if there exists a tilting object $T$ and a
nonempty open subset $U\subseteq\Fac(T,\alpha)$,
such that $\Coker\theta$ has property $P$ for all $\theta\in U$.
\end{dfn}

By Theorem \ref{t:change} this does not depend on the choice of objects
used to define $\Rep(T,\alpha)$.
In fact one does not need to consider all possible tilting objects,
but only one with $\Fac(T,\alpha)$ nonempty. Namely, we have
the following.

\begin{lem}
If the general object of dimension type $\alpha$ has property $P$,
and if $\Fac(T,\alpha)$ is nonempty, then there is a
nonempty open subset $U\subseteq \Fac(T,\alpha)$
such that $\Coker\theta$ has property $P$ for all $\theta\in U$.
\end{lem}

\begin{proof}
By definition this holds for some tilting object $T'$. Let $T''$ be a tilting
object with $T\oplus T'\in\Fac(T'')$.
There is a nonempty open subset $U'$ of $\Fac(T',\alpha)$ such that the corresponding
objects have property $P$. Since the property only depends on isomorphism classes
we may assume that $U'$ is $G(T',\alpha)$-stable.
Thus it corresponds to a nonempty $G(T'',\alpha)$-stable subset $U''$
which is open in $\Fac(T'',\alpha)_{T'}$, and hence in $\Fac(T'',\alpha)$.
Now $U''$ must meet the nonempty open set $\Fac(T'',\alpha)_{T}$,
and if their intersection is $V$, then $V$ corresponds to a nonempty
open subset $U$ of $\Fac(T,\alpha)$ with the required property.
\end{proof}

Similarly, if $Q$ is a property of pairs of objects $(X,Y)$,
depending only on the isomorphism classes, we
say that the \emph{general pair of objects $(X,Y)$,
with $X$ of dimension $\alpha$ and $Y$ of dimension $\beta$,
has property $Q$} if there exists a tilting object $T$ and a
nonempty open subset $U\subseteq\Fac(T,\alpha,\beta)$,
such that $(\Coker\theta,\Coker\phi)$ has property $Q$ for all $(\theta,\phi)\in U$.
The argument above shows the following.

\begin{lem}
If the general pair $(X,Y)$,
with $X$ of dimension $\alpha$ and $Y$ of dimension $\beta$,
have property $Q$, and if $\Fac(T,\alpha,\beta)$ is nonempty,
then there is a nonempty open subset $U\subseteq \Fac(T,\alpha,\beta)$
such that $(\Coker\theta,\Coker\phi)$ has property $Q$ for all $(\theta,\phi)\in U$.
\end{lem}

Given $\alpha,\beta\in K_0(\CC)_+$
we define $\hom(\alpha,\beta)$ to be the minimal possible value
of $\dim\Hom(X,Y)$ with $X$ of dimension type $\alpha$ and $Y$ of
dimension type $\beta$.
Similarly we define $\ext(\alpha,\beta)$ to be the
minimal possible value of $\dim\Ext^1(X,Y)$.
Clearly we have
\[
\hom(\alpha,\beta)-\ext(\alpha,\beta) = \langle \alpha,\beta\rangle.
\]
The upper semicontinuity of the dimension of Hom spaces, together
with the lemma above, show that the general pair of objects
$(X,Y)$, with $X$ of dimension $\alpha$ and $Y$ of dimension $\beta$,
has $\dim\Hom(X,Y)=\hom(\alpha,\beta)$ and $\dim\Ext^1(X,Y)=\ext(\alpha,\beta)$.

Let $\alpha=\beta_1+\beta_2+\dots+\beta_m$ with $\beta_i\in K_0(\CC)_+$.
If $T$ is a tilting object, then the set of $\theta\in\Fac(T,\alpha)$
such that $\Coker\theta$ is a direct sum of objects of dimensions
$\beta_1,\dots,\beta_m$ is the image of the map
\begin{equation}\label{e:dirsummap}
F : G(T,\alpha) \times \prod_{i=1}^m \Fac(T,\beta_i) \to \Fac(T,\alpha)
\end{equation}
defined by
\[
((f,g),(\theta_i)) \mapsto g \; \mathrm{diag}(\theta_1,\dots,\theta_m) \; f^{-1}.
\]
In particular this set is constructible. By considering possible
refinements of this decomposition (by a previous lemma one only needs to
consider finitely many such refinements), one sees that the set of $\theta$ with
$\Coker\theta$ a direct sum of \emph{indecomposable} objects of dimensions
$\beta_1,\dots,\beta_m$ is also constructible. The irreducibility of $\Fac(T,\alpha)$
then implies that any $\alpha\in K_0(\CC)_+$ has a
\emph{canonical decomposition} $\alpha=\beta_1+\beta_2+\dots+\beta_m$, unique up
to reordering, such that the general object of dimension
type $\alpha$ is a direct sum of indecomposable objects of
dimension types $\beta_1,\dots,\beta_m$.

\begin{thm}\label{t:dirext}
Let $\alpha=\beta_1+\beta_2+\dots+\beta_m$ with $\beta_i\in K_0(\CC)_+$.
The following are equivalent.
\begin{itemize}
\item[(i)]%
$\ext(\beta_i,\beta_j) = 0$ for $i\neq j$.
\item[(ii)]%
The general object of dimension $\alpha$ can be written as a direct
sum of objects of dimension types $\beta_1,\dots,\beta_m$.
\end{itemize}
\end{thm}

\begin{proof}
(i)$\Rightarrow$(ii).
Choose $T$ such that $\Rep(T,\beta_i)$ is nonempty for all $i$.
We may assume that $T'_\alpha = \bigoplus_{i=1}^m T'_{\beta_i}$
and $T''_\alpha = \bigoplus_{i=1}^m T''_{\beta_i}$.
By irreducibility we can find
\[
(\theta_1,\dots,\theta_m)\in\Fac(T,\beta_1)\times \dots\times\Fac(T,\beta_m)
\]
such that $\Ext^1(\Coker\theta_i,\Coker\theta_j)=0$ for all $i\neq j$.

Let $F$ be the map (\ref{e:dirsummap}), and consider the induced map $dF$ on
tangent spaces at the point $((1,1),(\theta_i))$.
It is the map
\[
\End(T'_\alpha) \oplus \End(T''_\alpha) \oplus \bigoplus_{i=1}^m \Hom(T'_{\beta_i},T''_{\beta_i})
\to
\Hom(T'_\alpha,T''_\alpha).
\]
given (in suggestive notation) by
\[
(df,dg,(d\theta_1,\dots,d\theta_m))
\mapsto
(dg_{ij}\theta_j - \theta_i df_{ij} + \delta_{ij} d\theta_i)_{ij}
\]
where we consider $df$ and $dg$ as matrices, with entries
$df_{ij}\in \Hom(T'_{\beta_j},T'_{\beta_i})$
and $dg_{ij}\in \Hom(T''_{\beta_j},T''_{\beta_i})$.

To compute the dimension of $\Ker dF$,
observe that if $i\neq j$, the $(i,j)$ component of this
matrix vanishes if and only if $(df_{ij},dg_{ij})\in \pi^{-1}(\theta_j,\theta_i)$,
and the vanishing of the $(i,i)$ component uniquely determines $d\theta_i$.
There is no restriction on the diagonal components $df_{ii}$ and $dg_{ii}$.
Thus
\[
\dim\Ker dF = \sum_{i\neq j}
\left( \dim\Hom(\Coker\theta_j,\Coker\theta_i)+\langle T''_{\beta_j},T'_{\beta_i}\rangle\right)
+
\]
\[
+ \sum_{i=1}^m \left( \langle T'_{\beta_i},T'_{\beta_i}\rangle + \langle T''_{\beta_i},T''_{\beta_i}\rangle\right).
\]
A dimension count shows that $dF$ is onto.
Thus $F$ is smooth at the point $(1,1,(\theta_i))$.
Thus $F$ is dominant.
Thus $\Fac(T,\alpha)$ contains a nonempty open subset $U$ such that
for all $\theta\in U$ one can write $\Coker\theta$ as a direct
sum of objects of dimensions $\beta_i$.
Thus (ii) holds.

\medskip
(ii)$\Rightarrow$(i).
Choose $T$ such that $\Fac(T,\alpha)$ contains a nonempty
open subset $U$ consisting of morphisms whose cokernel
is a direct sum of objects of dimension $\beta_i$.
Consider the nonempty open set $W$ of $\Fac(T,\alpha)$ consisting of the $\theta$
with $\dim\End(\Coker\theta)$ minimal.

By irreducibility the intersection $W\cap U$ must be nonempty,
so it contains a morphism $\theta$.
By assumption $\Coker\theta \cong X_1\oplus \dots \oplus X_m$
with $X_i$ of dimension type $\beta_i$.
Now if $i\neq j$ then $\Ext^1(X_i,X_j)=0$, for otherwise there is a nonsplit extension
$0\to X_j\to E\to X_i\to 0$, but then the long exact sequences
associated to this short exact sequence show that
\[
\dim\End(E\oplus \bigoplus_{k\neq i,j} X_k) < \dim\End(\bigoplus_{i=1}^m X_i).
\]
Since $E\oplus \bigoplus_{k\neq i,j} X_k$ is in $\Fac(T)$ and has
dimension type $\alpha$, this contradicts the definition of $W$.
\end{proof}

Note that since the argument of (i)$\Rightarrow$(ii) in the theorem shows
that the map (\ref{e:dirsummap}) is dominant, if the general objects
of dimension $\beta_i$ have properties $P_i$, then the general object
of dimension $\alpha$ can be written as a direct sum of objects of
dimension types $\beta_i$ satisfying $P_i$. It follows that
\begin{equation}\label{e:extsum}
\ext(\alpha,\gamma) = \sum_{i=1}^m \ext(\beta_i,\gamma)
\quad
\text{and}
\quad
\ext(\gamma,\alpha) = \sum_{i=1}^m \ext(\gamma,\beta_i)
\end{equation}
for all $\gamma\in K_0(\CC)_+$.

If $\alpha\in K_0(\CC)_+$ we say that $\alpha$ is \emph{generally indecomposable}
if the general object of dimension type $\alpha$ is indecomposable.
The following corollary is immediate.

\begin{cor}\label{c:genind}
A dimension type $\alpha$ is generally indecomposable
if and only if there is no non-trivial decomposition
$\alpha=\beta+\gamma$ with $0\neq \beta,\gamma\in K_0(\CC)_+$
and $\ext(\beta,\gamma) = \ext(\gamma,\beta) = 0$.
\end{cor}

\begin{proof}[Proof of Theorem \ref{t:can}]
If $\alpha=\beta_1+\beta_2+\dots$ is a decomposition
satisfying the conditions of Theorem \ref{t:dirext},
then some object has a decomposition
into direct summands, with dimension types given by the $\beta_i$,
and a decomposition into indecomposables,
with dimension types given by the canonical
decomposition $\alpha=\gamma_1+\gamma_2+\dots$.
By the Krull-Remak-Schmidt Theorem, we deduce that each $\beta_i$ is a sum of
a collection of the $\gamma_j$.
Now by Corollary \ref{c:genind} and (\ref{e:extsum}),
$\beta_i$ is generally indecomposable if and only
if this collection consists of only one $\gamma_j$.
\end{proof}

\section{Subobjects of general objects}\label{s:subobjects}
For $\alpha,\gamma\in K_0(\CC)_+$
we define $\gamma\hookrightarrow \alpha$ to
mean that the general object of dimension type $\alpha$
has a subobject of dimension type $\gamma$.

\begin{lem}
If $\gamma\hookrightarrow \alpha$ then there is a tilting object $S$ such that
$\Fac(S,\alpha)$ has a nonempty open subset $U$ such that for all $\theta\in U$,
$\Coker\theta$ has a subobject of dimension type $\gamma$ which belongs to $\Fac(S)$.
\end{lem}

\begin{proof}
Choose $T$ so that $\Fac(T,\alpha)$ is nonempty and let
$\beta=\alpha-\gamma$.
Let $L$ be the variety of
$((\theta,\phi),f,g)$ in $C(T,\alpha,\beta)$
with $(g \ \phi): T''_\alpha\oplus T'_\beta\to T''_\beta$
an epimorphism, or equivalently such that the
induced map $h$ on cokernels is an epimorphism.

By assumption image of the projection $p:L\to \Fac(T,\alpha)$ is
dense in $\Fac(T,\alpha)$.
Choose an irreducible locally closed subset $L'$ of $L$
which dominates $\Fac(T,\alpha)$.

Suppose $S$ is a tilting object with $\Fac(T)\subseteq \Fac(S)$.
We show that the set
\[
L_S = \{ ((\theta,\phi),f,g)\in L \mid \Ker h\in\Fac(S) \}
\]
is an open subset of $L$.
Now $\Ker h\in\Fac(S)$ if and only if $\Ext^1(S,\Ker h)=0$,
and thanks to the exact sequence
\[
0\to T'_\alpha\to \Ker(hp_\theta) \to \Ker h \to 0,
\]
it is equivalent that $\Ext^1(S,\Ker(hp_\theta))=0$.
Now if $((\theta,\phi),f,g)\in L$ then there is an
exact sequence
\[
0 \to \Ker(hp_\theta) \to T''_\alpha\oplus T'_\beta \to T''_\beta \to 0,
\]
so $\Ext^1(S,\Ker(hp))=0$ if and only if the induced map
\[
\Hom(S,(g\ \phi)):\Hom(S,T''_\alpha\oplus T'_\beta) \to \Hom(S,T''_\beta)
\]
is onto, and this is clearly an open condition.

Now $L = \bigcup_S L_S$, so $L'\cap L_S$ is an nonempty
open subset of $L'$ for some $S$, and hence it dominates $\Fac(T,\alpha)$.
Any element $\theta$ in the image has the property that $\Coker\theta$
has a subobject of dimension type $\gamma$ in $\Fac(S)$.
Now Theorem~\ref{t:change} enables one to pass from the open
subset of $\Fac(T,\alpha)$ to an open subset of $\Fac(S,\alpha)$.
\end{proof}

\begin{proof}[Proof of Theorem \ref{t:subext}]
Choose $T$ with $\Fac(T,\beta)$ and $\Fac(T,\gamma)$ nonempty.
We may assume that $T'_\alpha = T'_\beta\oplus T'_\gamma$
and $T''_\alpha = T''_\beta\oplus T''_\gamma$.

Suppose $Y$ is an object in $\Fac(T,\alpha)$ of dimension
type $\alpha$ and that $Y$ has a
subobject $X$ in $\Fac(T)$ of dimension type $\gamma$.
Then there is a commutative diagram
\[
\begin{CD}
& & 0 & & 0 & & 0 \\
& & @VVV @VVV @VVV \\
0 @>>> T'_\gamma @>\psi>> T''_\gamma @>p_\psi >> X @>>> 0 \\
& & @V(\begin{smallmatrix} 0 \\ 1 \end{smallmatrix})VV
@V(\begin{smallmatrix} 0 \\ 1 \end{smallmatrix})VV @VViV \\
0 @>>> T'_\beta\oplus T'_\gamma
@>\theta>> T''_\beta\oplus T''_\gamma @>(g\ ip_\psi)>> Y @>>> 0 \\
& & @V(1\ 0)VV @V(1\ 0)VV @VVpV \\
0 @>>> T'_\beta @>\phi>> T''_\beta @>p_\phi >> Y/X @>>> 0 \\
& & @VVV @VVV @VVV \\
& & 0 & & 0 & & 0
\end{CD}
\]
with exact rows and columns. Namely, $X$ and
the quotient $Y/X$ are both in $\Fac(T)$, so there
are $\phi\in\Fac(T,\beta)$ and $\psi\in\Fac(T,\gamma)$ giving
the top and bottom rows.
Since $X\in\Fac(T)$, the map
$\Hom(T''_\beta,Y)\to\Hom(T''_\beta,Y/X)$
is onto, so $p_\phi$ lifts to a morphism $g\in\Hom(T''_\beta,Y)$,
and this defines a morphism $(g\ ip_\psi):T''_\beta\oplus T''_\gamma\to Y$.
This gives the right hand half of the diagram.
Now the snake lemma gives an exact sequence of the kernels
\[
0\to T'_\gamma\to \Ker(g\ ip_\psi)\to T'_\beta\to 0,
\]
and this sequence must split, so it is equivalent to the split exact sequence
\[
0\to T'_\gamma\xrightarrow{(\begin{smallmatrix} 0 \\ 1 \end{smallmatrix})}
 T'_\beta\oplus T'_\gamma \xrightarrow{(1\ 0)} T'_\beta\to 0.
\]

Clearly the morphism $\theta$ in the diagram must be of shape
\[
\theta = \begin{pmatrix} \phi & 0 \\ \chi & \psi \end{pmatrix}
\]
for some $\chi\in\Hom(T'_\beta,T''_\gamma)$.
It follows that an object $Y$ in $\Fac(T)$ of dimension type $\alpha$
has a subobject in $\Fac(T)$ of dimension type $\gamma$
if and only if $Y$ is isomorphic to $\Coker\theta$ for some
$\theta$ of this shape, with $(\phi,\psi,\chi)$ in
\[
Z = \Fac(T,\beta)\times\Fac(T,\gamma)\times \Hom(T'_\beta,T''_\gamma).
\]

Consider the morphism
\[
F:G(T,\alpha)\times Z \to \Fac(T,\alpha),
\quad
((f,g),(\phi,\psi,\chi))\mapsto
g \begin{pmatrix} \phi & 0 \\ \chi & \psi \end{pmatrix} f^{-1}
\]
The induced map $dF$ on tangent spaces at $((1,1),(\phi,\psi,\chi))$
sends the element $((df,dg),(d\phi,d\psi,d\chi))$ to
\[
\begin{pmatrix}
d\phi + dg_{11}\phi + dg_{12}\chi - \phi df_{11}
&
dg_{12}\psi - \phi df_{12}
\\
d\chi + dg_{21}\phi + dg_{22}\chi - \chi df_{11} - \psi df_{21}
&
d\psi+dg_{22}\psi - \chi df_{12} - \psi df_{22}
\end{pmatrix}
\]
where we write $df$ and $dg$ as block matrices
\[
df = \begin{pmatrix} df_{11} & df_{12} \\ df_{21} & df_{22} \end{pmatrix},
\quad
dg = \begin{pmatrix} dg_{11} & dg_{12} \\ dg_{21} & dg_{22} \end{pmatrix}.
\]
To compute the kernel of $dF$, observe that $d\phi,d\psi,d\chi$ are
fixed by three of the entries of the matrix, and that the entries
of $df$ and $dg$ are unconstrained, except for the requirement that
$dg_{12}\psi = \phi df_{12}$, or equivalently that
$(df_{12},dg_{12})\in \pi^{-1}(\psi,\phi)$.
Thus
\[
\begin{split}
\dim \Ker dF = &\dim\Hom(\Coker\psi,\Coker\phi)+\langle T''_\gamma,T'_\beta\rangle + \\
&+ \langle T'_\alpha,T'_\alpha\rangle - \langle T'_\gamma,T'_\beta\rangle
+ \langle T''_\alpha,T''_\alpha\rangle - \langle T''_\gamma,T''_\beta\rangle.
\end{split}
\]
It follows that the cokernel of $dF$ has dimension
\[
\dim\Hom(\Coker\psi,\Coker\phi) +
\langle [T'_\gamma] - [T''_\gamma], [T''_\beta] - [T'_\beta]\rangle.
\]
Here the Euler form term is equal to $-\langle \gamma,\beta\rangle$,
so the cokernel of $dF$ has the same dimension as $\Ext^1(\Coker\psi,\Coker\phi)$.

Now if $\ext(\gamma,\beta) = 0$ then at some point $dF$ is onto, and hence
by Sard's lemma, $F$ is dominant, and so $\gamma\hookrightarrow\alpha$.
Conversely, if $\gamma\hookrightarrow\alpha$, then $F$ is dominant
so if the base field has characteristic zero then the induced map on
tangent spaces at some point of $G(T,\alpha)\times Z$ is surjective.
Clearly we may assume that this point has the form $((1,1),(\phi,\psi,\chi))$,
and then $\Ext^1(\Coker\psi,\Coker\phi)=0$, so that $\ext(\gamma,\beta) = 0$.
\end{proof}

\section{The general rank}\label{s:generalrank}
Given two objects $X$ and $Y$, we have a finite partition into disjoint locally closed
subsets
\[
\Hom(X,Y) = \bigcup_\gamma \Hom(X,Y)_\gamma,
\]
where $\Hom(X,Y)_\gamma$ denotes the set of morphisms of rank $\gamma\in K_0(\CC)$.
Thus there is a unique $\gamma$ such that $\Hom(X,Y)_\gamma$ is dense in $\Hom(X,Y)$.
We say that $\gamma$ is the \emph{general rank for $(X,Y)$}.

Given $\alpha,\beta\in K_0(\CC)_+$,
we say that the pair $(\alpha,\beta)$ has \emph{general rank $\gamma$}
if $\gamma$ is the general rank for the general pair $(X,Y)$, with $X$ of dimension
$\alpha$ and $Y$ of dimension $\beta$.
We show below that every pair $(\alpha,\beta)$ has a general rank, which is clearly unique.
We then use its existence to find a formula for $\hom(\alpha,\beta)$.

Let $T$ be a tilting object such that $\Fac(T,\alpha)$ and $\Fac(T,\beta)$
are both nonempty. We define
\[
\Fac(T,\alpha,\beta)' = \{ (\theta,\phi) \in \Fac(T,\alpha,\beta)
\mid
\dim\Hom(\Coker \theta,\Coker\phi)=\hom(\alpha,\beta)
\}
\]
By upper semicontinuity it is an open subset of $\Fac(T,\alpha,\beta)$,
and we have seen that it is nonempty. We also define
\[
C(T,\alpha,\beta)'=\pi^{-1}(\Fac(T,\alpha,\beta)').
\]

\begin{lem}
The map $C(T,\alpha,\beta)'\to \Fac(T,\alpha,\beta)'$ obtained by restricting $\pi$
is a vector bundle of rank $\hom(\alpha,\beta)+\langle T''_\alpha,T'_\beta\rangle$.
Thus $C(T,\alpha,\beta)'$ is irreducible.
\end{lem}

\begin{proof}
$C(T,\alpha,\beta)'$ is the kernel of the homomorphism of trivial vector bundles
\[
\Fac(T,\alpha,\beta)'\times \Hom(T'_\alpha,T'_\beta) \times\Hom(T''_\alpha,T''_\beta)
\to
\Fac(T,\alpha,\beta)'\times \Hom(T'_\alpha,T''_\beta),
\]
sending $((\theta,\phi),f,g)$ to $((\theta,\phi),g\theta-\phi f)$,
and its fibres have constant dimension by~(\ref{e:pifibses}).
Thus it is a sub-bundle.
\end{proof}

Given $\gamma$, let $C(T,\alpha,\beta)'_\gamma$ be the subset of $C(T,\alpha,\beta)'$
such that the induced morphism $h:\Coker\theta\to\Coker\phi$ has rank $\gamma$.

\begin{lem}
There are only finitely many $\gamma$ such that the sets $C(T,\alpha,\beta)'_\gamma$
are nonempty, and they are locally closed in $C(T,\alpha,\beta)'$.
\end{lem}

\begin{proof}
If $h$ is the induced morphism $\Coker\theta\to\Coker\phi$, then
we have an exact sequence $T''_\alpha\oplus T'_\beta \to T''_\beta \to \Coker h \to 0$.
Now the rank of $h$ is determined by the dimension type of $\Coker h$,
so by the rank of the morphism $T''_\alpha\oplus T'_\beta \to T''_\beta$.
Now the morphisms of any given rank form a locally closed subset
of $\Hom(T''_\alpha\oplus T'_\beta,T''_\beta)$, so they define locally
closed subsets of $C(\alpha,\beta)'$.
\end{proof}

It follows that if $\Fac(T,\alpha)$ and $\Fac(T,\beta)$
are both nonempty, then there is a unique $\gamma$ such that $C(T,\alpha,\beta)'_\gamma$
is a nonempty open subset of $C(T,\alpha,\beta)'$.
The next lemma shows that this is the general rank for
$(\alpha,\beta)$, so general ranks always exist.

\begin{lem}
If $C(T,\alpha,\beta)'_\gamma$ is a nonempty open subset of $C(T,\alpha,\beta)'$
then $\gamma$ is the general rank for $(\alpha,\beta)$.
\end{lem}

\begin{proof}
The map $C(T,\alpha,\beta)'\to \Fac(T,\alpha,\beta)'$ obtained
by restricting $\pi$ is a vector bundle, so
$U = \pi(C(T,\alpha,\beta)'_\gamma)$
is a nonempty open subset of $\Fac(T,\alpha,\beta)'$. Moreover,
if $(\theta,\phi)\in U$ then the set of $((\theta,\phi),f,g)\in C(T,\alpha,\beta)'_\gamma$
is a nonempty open subset of $\pi^{-1}(\theta,\phi)$.
Now for fixed $(\theta,\phi)\in U$, the map
\[
\pi^{-1}(\theta,\phi)\to \Hom(\Coker\theta,\Coker\phi)
\]
is a surjective linear map, hence a flat morphism of varieties. Thus the image of an open set is
open. Thus $\Hom(\Coker\theta,\Coker\phi)_\gamma$ is a nonempty
open subset of $\Hom(\Coker\theta,\Coker\phi)$. Thus
$\gamma$ is the general rank for $(\Coker\theta,\Coker\phi)$.
\end{proof}

\begin{lem}\label{l:genrankissubdim}
If $\gamma$ is the general rank for $(\alpha,\beta)$, then $\alpha-\gamma\hookrightarrow\alpha$ and
$\gamma\hookrightarrow \beta$.
\end{lem}

\begin{proof}
Clear.
\end{proof}

Now let $\gamma$ be the general rank for $(\alpha,\beta)$. Given a
tilting object $T$, we define $C(T,\alpha,\beta)''_\gamma$ to be
the set of $((\theta,\phi),f,g) \in C(T,\alpha,\beta)'_\gamma$
such that the induced morphism $h$ on cokernels has $\Ker h \in
\Fac(T)$.

\begin{lem}
$C(T,\alpha,\beta)''_\gamma$ is an open subset of $C(T,\alpha,\beta)'$.
\end{lem}

\begin{proof}
We know that $C(T,\alpha,\beta)'_\gamma$ is a nonempty open subset
of $C(T,\alpha,\beta)'$. There is a mapping
\[
C(T,\alpha,\beta)'_\gamma \to \Hom(\Hom(T,T''_\alpha\oplus T'_\beta),\Hom(T,T''_\beta)).
\]
sending $((\theta,\phi),f,g)$ to $\Hom(T,(g\ \phi))$.
If $h$ is the morphism on cokernels induced by $((\theta,\phi),f,g)$,
we show that
\[
\rank \Hom(T,(g\ \phi)) \le \langle T,[T''_\beta] - \beta + \gamma\rangle
\]
with equality if and only if $\Ker h\in \Fac(T)$,
from which the result follows.

Let $Z$ be the image of the morphism
$(g\ \phi):T''_\alpha\oplus T'_\beta\to T''_\beta$.
Breaking the sequence (\ref{e:gphiseq}) into two short exact sequences,
we obtain a long exact sequence
\[
0\to\Hom(T,\Ker(hp_\theta))\to \Hom(T,T''_\beta) \to \Hom(T,Z)\to \Ext^1(T,\Ker(hp_\theta)) \to 0,
\]
and
\[
0 \to\Hom(T,Z) \to \Hom(T,T''_\beta) \to \Hom(T,\Coker h)\to 0
\]
since $Z\in\Fac(T)$, so that $\Ext^1(T,Z)=0$.
Thus the rank of $\Hom(T,(g\ \phi))$ is the same as the rank
of the map
$\Hom(T,T''_\beta) \to \Hom(T,Z)$,
which is at most
\[
\dim\Hom(T,Z) = \langle T,Z\rangle
=\langle T,[T''_\beta] - \beta + \gamma\rangle.
\]
One has equality if and only if $\Ext^1(T,\Ker(hp_\theta))=0$.
Now sequence (\ref{e:kerhses}) gives a long exact sequence
\[
\dots \to \Ext^1(T,T'_\alpha) \to \Ext^1(T,\Ker(hp_\theta)) \to \Ext^1(T,\Ker h) \to \Ext^2(T,T'_\alpha)\to \dots
\]
and since $\Ext^1(T,T'_\alpha) = \Ext^2(T,T'_\alpha) = 0$ we
have $\Ext^1(T,\Ker(hp_\theta))=0$ if and only if $\Ext^1(T,\Ker h)=0$,
so if and only if $\Ker h\in\Fac(T)$.
\end{proof}

\begin{lem}
There is a tilting object $T$ such that $C(T,\alpha,\beta)''_\gamma$ is nonempty.
\end{lem}

\begin{proof}
Assuming just that $\Fac(S,\alpha)$ and $\Fac(S,\beta)$ are nonempty,
we know that $C(S,\alpha,\beta)'_\gamma$ is nonempty.
Let $((\theta,\phi),f,g)$ be an element of it, and let $h$ be the induced
morphism $\Coker\theta\to\Coker\phi$.
Clearly, if $T$ is a tilting object such that $\Fac(T)$ contains $S\oplus\Ker h$,
then $C(T,\alpha,\beta)''_\gamma$ is nonempty.
\end{proof}

Henceforth, let $T$ be a tilting object with $C(T,\alpha,\beta)''_\gamma$ nonempty.

\begin{lem}
If $((\theta,\phi),f,g)\in C(T,\alpha,\beta)''_\gamma$ and $h$ is the induced morphism on cokernels,
then the kernel of the epimorphism $T''_\alpha\to \Ima h$ is in $\add(T)$ and has dimension
type $[T''_\alpha] - \gamma$, and the kernel of the epimorphism $T''_\beta\to\Coker h$
is in $\add(T)$ and has dimension type $[T'_\beta] + \gamma$.
\end{lem}

\begin{proof}
The map $\Hom(T,T''_\alpha)\to \Hom(T,\Coker\theta)$ is onto since the next term in
the long exact sequence is $\Ext^1(T,T'_\alpha)=0$.
Also $\Hom(T,\Coker\theta)\to\Hom(T,\Ima h)$ is onto since the
next term in the long exact sequence is $\Ext^1(T,\Ker h)=0$.
Thus the map $\Hom(T,T''_\alpha)\to \Hom(T,\Ima h)$ is onto. Thus the
kernel of $T''_\alpha\to\Ima h$ is in $\add(T)$.

The map $\Hom(T''_\beta,\Coker\phi)\to \Hom(T''_\beta,\Coker h)$ is
onto, for the next term is $\Ext^1(T''_\beta,\Ima h)$. Thus the kernel
of $T''_\beta\to\Coker h$ is in $\add(T)$.
Its dimension type is
$[T''_\beta] - [\Coker h] = [T''_\beta] - [\Coker \phi] + \gamma = [T'_\beta]+\gamma$.
\end{proof}

Because of the lemma we know that there are objects $T'_*$ and $T''_*$ in $\add(T)$ with
\[
[T'_*] = [T''_\alpha] - \gamma,
\quad
[T''_*] = [T'_\beta] + \gamma.
\]

We define $V$ to be the set of tuples $(a,b,c,d,u,v)$ in
\begin{multline*}
\Hom(T'_*,T''_\alpha)\times  
\Hom(T'_*,T'_\beta)\times    
\Hom(T''_\alpha,T''_*) \times 
\Hom(T'_\beta,T''_*)\times \\ 
\Hom(T'_\alpha,T'_*)\times   
\Hom(T''_*,T''_\beta)        
\end{multline*}
such that $ca = db$ and such that $a,d,u,v$ are monomorphisms and
\[
(c\ d):T''_\alpha\oplus T'_\beta\to T''_*
\]
is an epimorphism.
This is clearly a locally closed subset, so a variety.

Observe that these conditions imply that $au$ and $vd$ are monomorphisms.
Also, the complex
\[
0\to T'_* \xrightarrow{(\begin{smallmatrix} -a \\ b\end{smallmatrix})}
T''_\alpha\oplus T'_\beta
\xrightarrow{(\begin{smallmatrix} c & d\end{smallmatrix})} T''_* \to 0
\]
is exact at the right hand term by definition, it is exact at the
left hand term since $a$ is a monomorphism, and then by dimensions
it is exact at the middle term.
Thus it is a (necessarily split) short exact sequence.
It follows that the map
\[
\Hom(T,(c\ d)):\Hom(T,T''_\alpha\oplus T'_\beta)\to \Hom(T,T''_*)
\]
is surjective.
Now since $v$ is a monomorphism,
one obtains a left exact sequence
\[
0\to T'_* \xrightarrow{(\begin{smallmatrix} -a \\ b\end{smallmatrix})}
T''_\alpha\oplus T'_\beta
\xrightarrow{(\begin{smallmatrix} vc & vd\end{smallmatrix})} T''_\beta.
\]

\begin{lem}
Suppose that $(a,b,c,d,u,v)\in V$ and let $\theta = au$, $\phi = vd$.
If $c_0,c_1,c_2,c_3$ are the induced maps on cokernels
\[
\begin{CD}
0 @>>> T'_\alpha @>u >> T'_* @>>> \Coker u @>>> 0 \\
& & @| @Va VV @Vc_0 VV \\
0 @>>> T'_\alpha @>\theta>> T''_\alpha @>p_\theta >> \Coker \theta @>>> 0 \\
& & @Vu VV @| @Vc_1 VV \\
0 @>>> T'_* @>a >> T''_\alpha @>>> \Coker a @>>> 0 \\
& & @Vb VV @Vvc VV @Vc_2 VV \\
0 @>>> T'_\beta @>\phi>> T''_\beta @>p_\phi >> \Coker \phi @>>> 0 \\
& & @Vd VV @| @Vc_3 VV \\
0 @>>> T''_* @>v >> T''_\beta @>>> \Coker v @>>> 0,
\end{CD}
\]
then the sequences
\[
0\to \Coker u\xrightarrow{c_0} \Coker \theta \xrightarrow{c_1} \Coker a\to 0
\]
and
\[
0\to \Coker a \xrightarrow{c_2} \Coker\phi \xrightarrow{c_3} \Coker v\to 0
\]
are exact.
\end{lem}

\begin{proof}
Although we are working in an abstract category $\CC$, as usual
we can verify this by diagram chasing.

$c_0$ is mono. (Diagram chase: if $w\in\Coker u$ is sent to 0, and $w$ comes from $t'_*$, then $a(t'_*) = \theta(t'_\alpha)$,
for some $t'_\alpha$. Then $a(t'_*-u(t'_\alpha))=0$. Thus since $a$ is mono, $t'_*=u(t'_\alpha)$. Thus $w=0$.

Exact at $\Coker\theta$.
If $c\in\Coker \theta$ is sent to 0 in $\Coker a$,
and $c$ comes from $t''_\alpha$, then $t''_\alpha = a(t'_*)$ for some
$t'_*$. But then the image of $t'_*$ in $\Coker u$ maps to $c$.

Clearly the maps $c_1$ and $c_3$ are epimorphisms.

$c_2$ is a monomorphism.
(Diagram chase: if $z\in \Coker a$ is sent to 0 by $c_2$,
and $z$ comes from $t''_\alpha$,
then the image of $t''_\alpha$ in $T''_\beta$ is also in
the image of $T'_\beta$, say from $t'_\beta$.
Then by the left exact sequence there is an element
$t'_*\in T'_*$ with $a(t'_*)=t''_\alpha$ and $b(t'_*)=t'_\beta$.
Then $t''_\alpha$ is sent to 0 in $\Coker a$.)

Exact at $\Coker\phi$.
(Diagram chase: if $z\in\Coker\phi$ is sent to 0,
and $z$ comes from $t''_\beta$, then $t''_\beta$
must come from an element $t''_*$. Then since $(-c\ d)$ is an epimorphism
there are $t''_\alpha$ and $t'_\beta$
with $t''_* = c(t''_\alpha)+d(t'_\beta)$.
Then $t''_\beta = v(t''_*) = vc(t''_\alpha)+\phi(t'_\beta)$,
so $z$ is the image of $vc(t''_\alpha)$.
Thus it is in the image of the map $\Coker a\to\Coker\phi$.)
\end{proof}

The lemma shows that there is a map
\[
f : V\to C(T,\alpha,\beta)''_\gamma,
\quad
f(a,b,c,d,u,v) = ((au,vd),bu,vc)
\]

\begin{lem}
The map $f$ is onto and every fibre has dimension
$\langle T''_*,T''_*\rangle + \langle T'_*,T'_*\rangle$.
\end{lem}

\begin{proof}
Suppose that $((\theta,\phi),f,g)\in C(T,\alpha,\beta)''_\gamma$,
and let $h$ be the induced map on cokernels.
The lemma shows that the kernel of the morphism $T''_\alpha\to\Ima h$ is isomorphic to $T'_*$,
and the kernel of the morphism $T''_\beta\to\Coker h$ is isomorphic to $T''_*$.
Thus we can find monomorphisms $a,v$ giving exact sequences
\[
0\to T'_* \xrightarrow{a} T''_\alpha \to \Ima h \to 0
\]
and
\[
0 \to T''_* \xrightarrow{v} T''_\beta \to \Coker h \to 0.
\]
In fact the set of possible $a,v$ is a variety of dimension $\langle T''_*,T''_*\rangle + \langle T'_*,T'_*\rangle$.
Then there are unique morphisms $b,d,u$ making the diagram
\[
\begin{CD}
0 @>>> T'_\alpha @>\theta>> T''_\alpha @>p_\theta >> \Coker \theta @>>> 0 \\
& & @Vu VV @| @VVV \\
0 @>>> T'_* @>a >> T''_\alpha @>>> \Ima h @>>> 0 \\
& & @Vb VV @VgVV @VVV \\
0 @>>> T'_\beta @>\phi>> T''_\beta @>p_\phi >> \Coker \phi @>>> 0 \\
& & @Vd VV @| @VVV \\
0 @>>> T''_* @>v >> T''_\beta @>>> \Coker h @>>> 0
\end{CD}
\]
commute. Clearly $u$ and $d$ are monomorphisms
since $\theta=au$ and $\phi=vd$. Now $f=bu$,
and the morphism from $T''_\alpha$ to $\Coker h$ is zero,
so there is a unique induced morphism $c\in\Hom(T''_\alpha,T''_*)$ with $g = vc$.
We have $vdb = \phi b = ga = vca$, and since $v$ is a
monomorphism this implies that $db=ca$.
Finally, consider the morphism
$(c\ d) : T''_\alpha\oplus T'_\beta \to T''_*$.
We show it is an epimorphism by diagram chasing.
Take $t''_*\in T''_*$. It gets sent by $v$ to $t''_\beta\in T''_\beta$.
This gets sent by $p_\phi$ to an element $z\in \Coker\phi$
whose image in $\Coker h$ is zero. Thus $z$ comes from
an element of $\Ima h$, so from an element $t''_\alpha\in T''_\alpha$.
Now $g(t''_\alpha)-t''_\beta$ is sent to 0 in $\Coker \phi$,
so $g(t''_\alpha)-t''_\beta=\phi(t'_\beta)$ for some $t'_\beta$.
Then $v(t''_*) = t''_\beta = g(t''_\alpha)-\phi(t'_\beta) = vc(t''_\alpha) - vd(t'_\beta)$.
Thus since $v$ is mono, $t''_* = c(t''_\alpha)-d(t'_\beta)$.
\end{proof}

\begin{lem}
Let $X$ be the set of $(c,d)$ with
$(c\ d):T''_\alpha\oplus T'_\beta\to T''_*$ an epimorphism
and $\Hom(T,(c\ d))$ surjective.
If $Y$ is the set of $(a,b,(c,d))$ with such that $ca = db$
and $(c,d)\in X$, then the projection $X\to Y$ is a sub-bundle
of the trivial vector bundle
\[
\Hom(T'_*,T''_\alpha)\times  
\Hom(T'_*,T'_\beta)\times   
X
\to X
\]
of rank $\langle T'_*,T'_* \rangle$.
Moreover $V$ is a nonempty open subset of
$Y\times\Hom(T'_\alpha,T'_*)\times\Hom(T''_*,T''_\beta)$.
\end{lem}

\begin{proof}
It is the kernel of the homomorphism of trivial bundles
\[
\Hom(T'_*,T''_\alpha)\times  
\Hom(T'_*,T'_\beta)\times   
X
\to
\Hom(T'_*,T''_*)\times X
\quad
(a,b,(c,d))\mapsto (ca-db,(c,d)).
\]
Since $(c\ d)$ is an epimorphism and $\Hom(T,(c\ d))$ is surjective,
one deduces that $\Ker (c\ d)\in\add(T)$. Thus, since it has dimension
\[
[T''_\alpha] + [T'_\beta] - [T''_*] = [T'_*],
\]
we have $\Ker(c\ d)\cong T'_*$. Now the fibre of $Y$ over $(c\ d)$
is identified with the space $\Hom(T'_*,\Ker(c\ d)) \cong \End(T'_*)$,
so it has constant dimension $\langle T'_*,T'_* \rangle$.
\end{proof}

We can now prove an analogue of Schofield's
formula \cite[Theorem 5.2]{S} in our setting.

\begin{thm}\label{t:genrankformula}
If $\gamma$ is the general rank for $\alpha,\beta$, then
\[
\ext(\alpha,\beta) = -\langle \alpha-\gamma,\beta-\gamma\rangle.
\]
\end{thm}

\begin{proof}
$V$ is nonempty since $f$ is onto. We compute its dimension
in two ways.

By the lemma above $V$ is an open subset of
a known irreducible variety. Thus $V$ is irreducible of dimension
\[
\langle T''_\alpha\oplus T'_\beta,T''_*\rangle
+ \langle T'_*,T'_*\rangle
+\langle T'_\alpha,T'_* \rangle
+\langle T''_*,T''_\beta\rangle.
\]
With $\alpha = T''_\alpha - T'_\alpha$, $\beta= T''_\beta - T'_\beta$
and $\gamma = T''_\alpha- T'_* = T''_* - T'_\beta$
this becomes
\[
\langle T''_\alpha+T''_\beta-\beta,T''_\beta+\gamma-\beta\rangle
+ \langle T''_\alpha-\gamma,T''_\alpha-\gamma\rangle
+\langle T''_\alpha-\alpha,T''_\alpha-\gamma \rangle
+\langle T''_\beta+\gamma-\beta,T''_\beta\rangle.
\]

On the other hand,
$\Fac(T,\alpha,\beta)'$ is open in $\Fac(T,\alpha,\beta)$,
so it is irreducible of dimension
$d = \langle T'_\alpha,T''_\alpha\rangle+\langle T'_\beta,T''_\beta\rangle$.
Now $C(T,\alpha,\beta)'$ is a vector bundle over $\Fac(T,\alpha,\beta)'$
of rank $r = \hom(\alpha,\beta)+\langle T''_\alpha,T'_\beta\rangle$, so
it is irreducible of dimension $d+r$.
Since $C(T,\alpha,\beta)''_\gamma$ is an open subset of $C(T,\alpha,\beta)'$,
so it too is irreducible of dimension $d+r$.
Then the morphism $f:V\to C(T,\alpha,\beta)''_\gamma$ is onto and has
fibres of dimension $k = \langle T''_*,T''_*\rangle + \langle T'_*,T'_*\rangle$,
so $V$ has dimension $d+r+k$.

Equating these two expressions gives the required result.
\end{proof}

\begin{proof}[Proof of Theorem \ref{t:extform}]
If $\delta\hookrightarrow\alpha$ and $\eta\hookrightarrow\beta$ then there are
$X$ and $Y$ of dimension types $\alpha$ and $\beta$ respectively, with
$\dim\Ext^1(X,Y)=\ext(\alpha,\beta)$ and such that $X$ and $Y$ have
subobjects $D$ and $H$ of dimensions $\delta$ and $\eta$ respectively.
Since $\CC$ is hereditary, the maps $\Ext^1(X,Y)\to \Ext^1(D,Y)$ and
$\Ext^1(D,Y)\to \Ext^1(D,Y/H)$ are onto, so
\begin{align*}
\ext(\alpha,\beta) &= \dim\Ext^1(X,Y) \ge \dim\Ext^1(D,Y/H)
\\
&= \dim\Hom(D,Y/H)-\langle \delta,\beta-\eta\rangle
\ge -\langle \delta,\beta-\eta\rangle.
\end{align*}
The result now follows from Theorem~\ref{t:genrankformula} and Lemma~\ref{l:genrankissubdim}.
\end{proof}

\begin{proof}[Proof of Corollary \ref{c:schur}]
(4) $\Rightarrow$ (3).
If there is an expression $\alpha=\beta+\gamma$ with $\beta\hookrightarrow\alpha$ and
$\gamma\hookrightarrow\alpha$, apply (4) to both $\beta$ and $\gamma$
to get a contradiction.

(3) $\Rightarrow$ (2).
Choose $T$ such that $\Fac(T,\alpha)$ is nonempty.
Suppose that for all $\theta\in\Fac(T,\alpha)$, the
object $\Coker\theta$ has nontrivial endomorphism algebra.
Then since $K$ is algebraically closed, it has an endomorphism
of rank $\delta\neq \alpha,0$.

Let $W$ be the subset of $C(T,\alpha,\alpha)'$ consisting of quadruples of the form
$(\theta,\theta,f,g)$.

Let $f:W\to\Fac(T,\alpha)$ be the projection.

By assumption $\bigcup_{\delta\neq\alpha,0} (W\cap C(T,\alpha,\alpha)'_\delta)$
dominates $\Fac(T,\alpha)$. Thus some $W\cap C(T,\alpha,\alpha)'_\delta$
dominates $\Fac(T,\alpha)$. Now every $\theta$ in the image has the property
that $\Coker\theta$ has an endomorphism of rank $\delta$. Thus it has
subobjects of dimensions $\delta$ and $\alpha-\delta$. This contradicts (3).

(2) $\Rightarrow$ (1) is clear.

(1) $\Rightarrow$ (3).
If there were such an expression, then (using characteristic 0)
we get that $\ext(\beta,\gamma) = \ext(\gamma,\beta) = 0$,
contradicting that the canonical decomposition of $\alpha$ is $\alpha$ itself.

(2) $\Rightarrow$ (4).
Suppose that $\beta\hookrightarrow\alpha$.
By irreducibility we can find an object $X$ of dimension $\alpha$ with trivial
endomorphism algebra, which has a subobject $Y$ of dimension $\beta$.
Now $\Hom(X/Y,Y)=0$ since $X$ has trivial endomorphism algebra.
Also, there is a nonsplit extension $0\to Y\to X\to X/Y\to 0$,
so that $\Ext^1(X/Y,Y)\neq 0$. Thus $\langle\alpha-\beta,\beta\rangle<0$.
On the other hand, since characteristic zero and $\beta\hookrightarrow\alpha$
we have $\ext(\beta,\alpha-\beta)=0$, and hence $\langle\beta,\alpha-\beta\rangle\ge 0$.
\end{proof}

\section{Weighted projective lines}\label{s:wtd}
We use the notation and terminology of Geigle and Lenzing \cite{GL}
throughout this section.
Let $\mathbb{X}$ be the weighted projective line associated to
a weight sequence $\pp=(p_0,\dots,p_n)$ and distinct points
$\blambda=(\lambda_0,\dots,\lambda_n)$ in $\PP^1(K)$.
Let $\CC$ be the category of coherent sheaves on $\mathbb{X}$.
It is shown in \cite{GL} that $\CC$ has property (H).

Recall that there is a rank one additive abelian group $\LL(\pp)$
with generators $\vec x_0,\vec x_1,\dots,\vec x_n$ and relations
$p_0 \vec x_0 = p_1 \vec x_1 = \dots = p_n\vec x_n = \vec c$, say,
and that the line bundles on $\XX$ are classified as $\OO(\vec{x})$
with $\vec{x}\in\LL(\pp)$. Moreover there is a partial order
on $\LL(\pp)$ given by $\vec x \le \vec y$ if and only if
\[
\vec y - \vec x = \sum_{i=0}^n k_i \vec x_i
\]
for some $k_i\ge 0$, and
\[
\Hom(\OO(\vec{x}),\OO(\vec{y}))\neq 0 \Longleftrightarrow \vec{x}\le\vec{y}.
\]

For each ordinary point of $\PP^1(K)$ there is a unique simple
torsion sheaf $S$ with support at this point, and there
is an exact sequence
\[
0\to \OO(0)\to \OO(\vec c)\to S\to 0,
\]
while for each exceptional point $\lambda_i$ there
are simple sheaves $S_{ij}$ ($0\le j<p_i$) concentrated at the point,
and exact sequences
\[
0\to \OO((j-1)\vec x_i)\to \OO(j\vec x_i)\to S_{ij}\to 0.
\]
(We number the simple sheaves following Schiffmann \cite[\S 2.4]{Schiffmann} rather than \cite{GL}.)

It is shown in \cite[Proposition 4.1]{GL} that
\[
T_0 = \bigoplus_{0\le \vec{x} \le \vec{c}} \OO(\vec{x})
= \OO(0) \oplus \left( \bigoplus_{i=0}^n \bigoplus_{j=1}^{p_i-1} \OO(j \vec x_i) \right) \oplus \OO(\vec c)
\]
is a tilting sheaf.

Clearly the shifts $T_m = T_0(m\vec{c})$ ($m\in\Z$) are tilting
sheaves. Now
\[
\Ext^1(T_{m-1},T_m) \cong \Hom(T_m,T_{m-1}(\vec{\omega}))^* \cong \Hom(T_0,T_0(\vec{\omega}-\vec{c}))^*
\]
by Serre duality.
If this were nonzero, then by the above one would have $\vec{x}\le \vec{\omega}-\vec{c}+\vec{y}$ for
some $0\le \vec x,\vec y\le \vec c$, and hence $0\le \vec\omega$, which is not true.
Thus $\Ext^1(T_{m-1},T_m)=0$, so
\[
\dots \supset \Fac(T_{-1}) \supset \Fac(T_0) \supset \Fac(T_1) \supset \dots
\]
Also, for any sheaf $X$ one has $\Ext^1(T_m,X) \cong \Ext^1(T_0,X(-m))$, and this
vanishes for $m\ll 0$, so that $X\in \Fac(T_m)$. Thus $\CC$ has property (T).
It remains to prove that $\CC$ has property (R).

\begin{lem}
There is a group homomorphism, $\partial:K_0(\CC)\to \LL(\pp)$
such that $\partial([\OO(\vec x)]) = \vec x$
for all $\vec x\in\LL(\pp)$,
$\partial([S_{ij}])=\vec x_i$ for all $i,j$,
and $\partial([S])=\vec c$ if $S$ is a simple torsion sheaf
concentrated at a non-exceptional point.
\end{lem}

We call $\partial$ the \emph{weighted degree}.

\begin{proof}
The dimension types of the indecomposable summands of $T_0$
freely generate $K_0(\CC)$, so we can define
$\partial$ to be the unique homomorphism with $\partial([\OO(\vec x)]) = \vec x$
for $0\le \vec x\le \vec c$. The exact sequences above show that $\partial$
has the right effect on the dimension types of the simple torsion sheaves.
Now for any $\vec x$ there is an exact sequence
\[
0\to \OO(\vec x)\to \OO(\vec x+\vec x_i)\to S_{i1}(\vec x)\to 0
\]
and $S_{i1}(\vec x)\cong S_{ij}$ for some $j$.
Thus $\partial([\OO(\vec x+\vec x_i)]) = \partial([\OO(\vec x)])+\vec x_i$.
It follows that $\partial([\OO(\vec x)]) = \vec x$ for all $\vec x$.
\end{proof}

Note that the weighted degree does not
respect the partial orderings on $K_0(\CC)$ and $\LL(\pp)$.

Clearly the weighted degree of any torsion sheaf is $\ge 0$.

\begin{lem}
If $T'\in\add(T_m)$ and $Y$ is torsion of weighted degree $\le k\vec c$, then the
kernel of any morphism $\theta:T'\to Y$ is in $\Fac(T_{m-k})$.
\end{lem}

\begin{proof}
We prove this by induction on the rank of $T'$.
If $T'$ is a line bundle, say $\OO(\vec x)$, then the image of
$\theta$ is torsion of weighted degree $\vec z\le k\vec c$, say, and
$\Ker\theta\cong \OO(\vec x-\vec z) \in \Fac(T_{m-k})$, as required.
Otherwise choose a proper direct summand $T''$ of $T'$. Then there is an
exact sequence
\[
0\to \Ker(\theta|_{T''}) \to \Ker\theta \to \Ker\phi\to 0
\]
where $\phi$ is the induced morphism $T'/T'' \to Y/\theta(T'')$.
By induction $\Ker(\theta|_{T''})$ and $\Ker\phi$ are in
$\Fac(T_{m-k})$, and hence so is $\Ker\theta$.
\end{proof}

\begin{lem}
Let $m\in\Z$ and $q\ge 0$.
If $Y$ is a line bundle direct summand of $T_{m+q}$ then any morphism
$\theta:T'\to Y$ with $T'\in\add(T_m)$ has kernel in $\Fac(T_{m-q-1})$.
\end{lem}

\begin{proof}
If $\theta=0$ this is clear, so suppose $\theta\neq 0$.
Then the restriction of $\theta$ to some line bundle
direct summand $T'_0$ is nonzero. Let $T'_1$ be the complementary
summand.
Let $T'_0 \cong \OO(\vec x)$ and $Y \cong \OO(\vec y)$.
The exact sequence
\[
0\to T'_0 \to Y\to Y/\theta(T'_0)\to 0,
\]
shows that $Y/\theta(T'_0)$ is torsion of weighted degree
\[
\vec y - \vec x \le (m+q+1)\vec c - m\vec c = (q+1)\vec c
\]
Now there is an exact sequence
\[
0\to\Ker(\theta|_{T'_0}) \to \Ker\theta\to \Ker\phi\to 0
\]
where $\phi$ is the morphism $T'_1\to Y/\theta(T'_0)$ induced by $\theta$.
However $\Ker(\theta|_{T'_0})=0$, so $\Ker\theta\cong \Ker\phi\in \Fac(T_{m-q-1})$
by the previous lemma.
\end{proof}

\begin{thm}
For any sheaves $X$ and $Y$ there is a tilting sheaf $T_s$ such
that $\Ker\theta\in\Fac(T_s)$ for all $\theta\in\Hom(X,Y)$.
\end{thm}

\begin{proof}
Let us say that a sheaf $Y$ is \emph{good} if for any
$m$ there is a tilting sheaf $T_s$ such that
for any morphism $\theta:X\to Y$ with $X\in\Fac(T_m)$,
one has $\Ker\theta\in\Fac(T)$.
Clearly it suffices to prove that every sheaf is good.

If $\theta\in\Hom(X,Y)$ and $X\in\Fac(T_m)$ then there is an
epimorphism $\psi:T'\to X$ with $T'\in\add(T_m)$, and
this induces an epimorphism $\Ker (\theta\psi) \to \Ker\theta$.
Thus the last two lemmas show that torsion sheaves
and line bundles are good.

Clearly any subsheaf of a good sheaf is good.
Since any vector bundle embeds in a direct sum of line bundles,
and any coherent sheaf is a direct sum of a vector bundle
and a torsion sheaf, it suffices to show that
if $Y'$ and $Y''$ are good, then so is $Y'\oplus Y''$.

Given $m$, since $Y'$ is good, there is a tilting sheaf
$T_r$ such that the kernel of any morphism from a sheaf
in $\Fac(T_m)$ to $Y'$ is in $\Fac(T_r)$.
Since $Y''$ is good, there is a tilting sheaf $T_s$ such
that kernel of any morphism from a sheaf in $\Fac(T_r)$
to $Y'$ is in $\Fac(T_s)$. Now if $X\in\Fac(T_m)$ and
\[
\theta = \begin{pmatrix} \theta' \\ \theta'' \end{pmatrix} :X\to Y'\oplus Y''
\]
then $\Ker \theta'\in\Fac(T_r)$, and then since there is an exact sequence
\[
0\to \Ker\theta \to \Ker\theta'\xrightarrow{\phi} Y'',
\]
where $\phi$ is the restriction of $\theta''$, we have $\Ker\theta\in\Fac(T_s)$.
\end{proof}

\begin{thm}
$\CC$ satisfies condition \textup{(R)}.
\end{thm}

\begin{proof}
Given $X,Y$, let $T$ be a tilting sheaf with $Y\in\Fac(T)$ and
$\Ker\theta\in\Fac(T)$ for all $\theta\in\Hom(X,Y)$.
Let $T_i$ be the non-isomorphic indecomposable summands of $T$.

Let $\theta\in\Hom(X,Y)$. The natural map $\Hom(T_i,\Ima\theta)\to\Hom(T_i,Y)$ is 1-1,
and since $\Ker\theta\in\Fac(T)$ the map $\Hom(T_i,X)\to\Hom(T_i,\Ima\theta)$ is
onto. It follows that
\[
\rank \Hom(T_i,\theta) = \dim \Hom(T_i,\Ima \theta) = \langle T_i, \Ima \theta\rangle.
\]
Since the elements $[T_i]$ generate $K_0(\CC)$, and the Euler form is
non-degenerate, the numbers $\langle T_i, \Ima \theta\rangle$
determine the rank of $\theta$.
Since there are only finitely many possible ranks of $\Hom(T_i,\theta)$,
there are only finitely many possible ranks of morphisms $\theta$.

Now the map
\[
\Hom(X,Y) \to \Hom_K(\Hom(T_i,X),\Hom(T_i,Y))
\]
sending $\theta$ to $\Hom(T_i,\theta)$ is linear, so it is a morphism
of varieties, and hence the set of $\theta$ for which $\Hom(T_i,\theta)$
has any given rank is locally closed. It follows that $\Hom(X,Y)_\gamma$
is locally closed in $\Hom(X,Y)$.

The set of monomorphisms is given by the $\theta$ with $\Hom(T_i,\theta)$
having rank $\Hom(T_i,X)$, and clearly this is an open condition.
The set of epimorphisms is given by the $\theta$ with $\Hom(T_i,\theta)$
having rank $r = \langle T_i, Y\rangle$, and since we have
assumed that $Y\in\Fac(T)$, we have $r=\dim\Hom(T_i,Y)$,
so the set of such $\theta$ is also open.
\end{proof}

\section{An algorithm for weighted projective lines}\label{s:algorithm}
Following Schiffmann \cite{Schiffmann}, it is useful
to identify $K_0(\CC)$ with the root lattice of a loop algebra of a Kac-Moody Lie algebra,
and take as $\Z$-basis the elements $\alpha_* = [\OO(0)]$,
$\alpha_{ij} = [S_{ij}]$ for $0\le i\le n$ and $1\le j\le p_i-1$, and
$\delta = [\OO(\vec c)]-[\OO(0)]$.
We write $\alpha\in K_0(\CC)$ as
\begin{equation}\label{e:alphaform}
\alpha = m_* \alpha_* + \sum_{i=0}^n\sum_{j=1}^{p_i-1} m_{ij} \alpha_{ij} + d\delta
\end{equation}
with $m_*,m_{ij},d\in\Z$. If $\beta$ is given by the corresponding formula with
coefficients $m'_*,m'_{ij},d'$, the Euler form is given by
\[
\langle \alpha,\beta\rangle = m_*m'_*
+ m_* d' - dm'_*
+ \sum_{i=0}^n \sum_{j=1}^{p_i-1} m_{ij}(m'_{ij} - m'_{i,j-1})
\]
with the convention that $m'_{i0} = m'_*$, see \cite[Lemma 6.1]{Schiffmann}.

We write $K_0(\CC)_{\mathrm{tors}}$ for the set of dimension types of torsion sheaves.

\begin{lem}
The set $K_0(\CC)_{\mathrm{tors}}$ consists of the $\alpha$ with $m_*=0$
and
\begin{equation}\label{e:dineq}
d \ge \sum_{i=0}^n \max \{ 0,-m_{ij} \mid 1\le j\le p_i-1\}.
\end{equation}
\end{lem}

\begin{proof}
$K_0(\CC)_{\mathrm{tors}}$ consists of the
non-negative linear combinations of $\delta$, $\alpha_{ij}$
and $[S_{i0}] = \delta - \sum_{j=1}^{p_i-1} \alpha_{ij}$,
say
\begin{equation}\label{e:torsionform}
k\delta
+ \sum_{i=0}^n \sum_{j=1}^{p_i-1} k_{ij}\alpha_{ij}
+ \sum_{i=0}^n k_i \biggl(\delta - \sum_{j=1}^{p_i-1} \alpha_{ij}\biggr)
= d\delta + \sum_{i=0}^n \sum_{j=1}^{p_i-1} m_{ij} \alpha_{ij}
\end{equation}
where $d = k+\sum_i k_i$ and $m_{ij} = k_{ij}-k_i$.
Clearly $k_i = k_{ij} - m_{ij} \ge -m_{ij}$, so $k_i \ge \max\{0,-m_{ij}\}$,
and the inequality for $d$ follows. Conversely if the inequality for $d$
holds, then $\alpha$ arises on taking $k_i = \max\{0,-m_{ij} \}$,
$k = d-\sum_i k_i$ and $k_{ij} = m_{ij}+k_i$.
\end{proof}

\begin{lem}
The set $K_0(\CC)_+$ consists of the $\alpha$ with $m_* > 0$ or with $m_*=0$
and $d$ satisfying \textup{(\ref{e:dineq})}.
\end{lem}

\begin{proof}
If $\alpha=[X]$ then $m_*$ is the rank of $X$, so it is non-negative.
If $m_*>0$, then for $N\gg 0$ there is a torsion sheaf of dimension type
$\alpha - m_*\alpha_* + Nm_* \delta$, and its direct sum with $m_*$
copies of $\OO(-N\vec c)$ then has dimension type $\alpha$.
\end{proof}

The vector bundles in $\CC$ can be identified with vector bundles
$E$ on $\PP^1$ with a quasi-parabolic structure of type $\pp$, that is,
a flag of subspaces
\[
E_{\lambda_i} \supseteq E_{i1} \supseteq E_{i2} \supseteq \dots \supseteq E_{i,p_i-1},
\]
of the fibre at each exceptional point $\lambda_i$, see \cite[\S 4.2]{L}.
This identification is not unique, but it can be done in such a
way that if $\vec x\in \LL(\pp)$ is written in the form $\vec x = \ell \vec c + \sum_{i=0}^n \ell_i \vec x_i$,
with $0\le \ell_i<p_i$, then $\OO(\vec x)$ corresponds to the bundle $E=\OO(\ell)$ on $\PP^1$
with $E_{ij}$ one dimensional for $j\le\ell_i$ and zero for $j>\ell_i$.
It follows that the dimension type can be written as
\[
(\rank E) \alpha_* + \sum_{i=0}^n \sum_{j=1}^{p_i-1} (\dim E_{ij}) \alpha_{ij} + (\deg E)\delta.
\]
Thus the possible dimension types of vector bundles
in $\CC$ are the $\alpha$ with
\[
m_* \ge m_{i1} \ge m_{i2} \ge \dots \ge m_{i,p_i-1} \ge 0
\]
for all $i$, and $d=0$ if $m_*=0$. We denote the set of these $\alpha$
by $K_0(\CC)_{\mathrm{vb}}$.
If $\alpha$ belongs to this set then the general sheaf of dimension
type $\alpha$ is a vector bundle, for, if $X$ is a vector bundle
of dimension $\alpha$,  there is some $h\gg 0$ with $\Hom(T_h,X)=0$,
and then $\Hom(T_h,Y)=0$ for the general sheaf $Y$ of dimension
type $\alpha$, so that $Y$ is a vector bundle.

\begin{thm}\label{t:torsdecomp}
If $\alpha\in K_0(\CC)_+$, then the general sheaf of dimension type
$\alpha$ is the direct sum of a vector bundle of dimension type $\beta$
and a torsion sheaf of dimension type $\gamma = \alpha-\beta$, where
\begin{align*}
\beta &= m_* \alpha_* + \sum_{i=0}^n\sum_{j=1}^{p_i-1} m'_{ij} \alpha_{ij} + (d-\sum_{i=0}^n k_i)\delta,\\
k_i &= \max \{ 0,-m_{ij} \mid 1\le j\le p_i-1\},\\
m'_{ij} &= \min\{m_*,m_{i1}+k_i,m_{i2}+k_i,\dots,m_{ij}+k_i\}
\end{align*}
\end{thm}

We call $\alpha=\beta+\gamma$ the \emph{torsion} decomposition of $\alpha$.

\begin{proof}
Clearly $\beta\in K_0(\CC)_{\mathrm{vb}}$, and setting $k_{ij} = m_{ij}+k_{i} - m'_{ij} \ge 0$,
we have
\[
\gamma = \sum_{i=0}^n k_i \biggl(\delta - \sum_{j=1}^{p_i-1} \alpha_{ij}\biggr)
+ \sum_{i=0}^n \sum_{j=1}^{p_i-1} k_{ij} \alpha_{ij} \in K_0(\CC)_{\mathrm{tors}}.
\]
Now
\begin{align*}
\langle \gamma,\beta\rangle
&= -(\sum_{i=0}^n k_i) m_*
+ \sum_{i=0}^n \sum_{j=1}^{p_i-1} (k_{ij} - k_i)(m'_{ij} - m'_{i,j-1})
\\
&= \sum_{i=0}^n \sum_{j=1}^{p_i-1} k_{ij} (m'_{ij} - m'_{i,j-1}) - \sum_{i=0}^n k_i m'_{i,p_i-1}.
\end{align*}
We show that all terms in this expression are zero.
The definition of $m'_{ij}$ immediately gives $m'_{ij} = \min\{ m'_{i,j-1},m_{ij}+k_i\}$.
Now if $k_{ij}>0$, then $m'_{ij} < m_{ij}+k_i$, so $m'_{ij} = m'_{i,j-1}$, and hence
the terms in the first sum are zero.
If $k_i>0$, then there is some $m_{ij}<0$ with $k_i=-m_{ij}$. Then $m_{ij}+k_i=0$,
which implies that $m'_{i,p_i-1}=0$, and so the terms in the second sum are zero.

Thus $\langle \gamma,\beta\rangle = 0$.
Since there is a torsion sheaf of dimension $\gamma$ and a vector
bundle of dimension $\beta$, we have $\ext(\beta,\gamma)=0$
and $\hom(\gamma,\beta)=0$. Thus also $\ext(\gamma,\beta)=0$.
The assertion now follows from Theorem \ref{t:dirext} and the remark following it.
\end{proof}

\begin{lem}\label{l:underbun}
If $\alpha \in K_0(\CC)_{\mathrm{vb}}$,
then the general sheaf of dimension type $\alpha$
is a vector bundle whose underlying vector bundle on $\PP^1$ is isomorphic to
\[
\OO(t)^{\oplus a}\oplus \OO(t+1)^{\oplus b},
\]
where $a,b\ge 0$ satisfy $a+b=m_*$
and $ta+(t+1)b=d$.
\end{lem}

Observe that $t = \Floor(d/m_*)$, $a = (t+1)m_*-d$ and $b = d-t m_*$.

\begin{proof}
Let $Y$ (respectively $Z$) be the line bundle with
underlying vector bundle is $\OO(t-1)$ (respectively $\OO(t+1)$) and all
subspaces in the quasi-parabolic structure 1-dimensional (respectively zero).
There is a vector bundle $X$ of dimension type $\alpha$ whose underlying
vector bundle is isomorphic to $\OO(t)^{\oplus a}\oplus \OO(t+1)^{\oplus b}$.
Then we must have $\Hom(X,Y) = \Hom(Z,X)=0$, and so by semicontinuity
$\Hom(X',Y) = \Hom(Z,X')=0$ for the general representation $X'$ of dimension type $\alpha$.
Now any homomorphism of the underlying vector bundles on $\PP^1$
necessarily respects these quasi-parabolic structures, so there
are no homomorphisms from the underlying vector bundle of $X'$
to $\OO(t-1)$ or from $\OO(t+1)$.
The result thus follows from the Grothendieck-Birkhoff Theorem.
\end{proof}

\begin{proof}[Proof of Theorem \ref{t:recalg}]
By Theorem~\ref{t:extform}, we have
\begin{equation}\label{e:rankformeta}
\ext(\alpha,\beta) = \min \{ - \langle \alpha-\eta,\beta-\eta\rangle
\mid \alpha-\eta \hookrightarrow \alpha, \eta\hookrightarrow \beta \}
\end{equation}
which is useful provided that we can determine all $\eta$ with
$\alpha-\eta \hookrightarrow \alpha$ and $\eta\hookrightarrow \beta$.
Now by Theorem~\ref{t:subext}, we have
$\alpha-\eta \hookrightarrow \alpha$ if and only if $\ext(\alpha-\eta,\eta)=0$ and
$\eta\hookrightarrow \beta$ if and only if $\ext(\eta,\beta-\eta)=0$.
This would give a recursive algorithm, except that in this
simple form it doesn't terminate.

If $\alpha$ and $\beta$ are torsion dimension types, however, it does.
We use the degree in the sense of \cite[Proposition 2.8]{GL},
\[
\deg(\alpha) = k + \sum_{i=0}^n \biggl( k_i + \sum_{j=1}^{p_i-1} k_{ij} \biggr) \frac{p}{p_i},
\]
where $p = \lcm\{p_0,p_1,\dots,p_n\}$ and $\alpha$ is written in the form (\ref{e:torsionform}).
Observe that $\deg(\alpha)$ is a positive integer if $\alpha$ is nonzero.
There are only finitely many possible $\eta$ in (\ref{e:rankformeta})
as they must have $\deg(\eta) \le \min\{ \deg(\alpha),\deg(\beta)\}$,
and the recursion terminates by consideration of $\deg(\alpha+\beta)$.

In general we use Theorem~\ref{t:torsdecomp} to reduce to the case
when $\alpha,\beta\in K_0(\CC)_{\mathrm{tors}} \cup K_0(\CC)_{\mathrm{vb}}$,
and consider the rank of $\alpha+\beta$.

If $\alpha$ is a vector bundle dimension type and $\beta$ is torsion, then
$\ext(\alpha,\beta)=0$.
If $\alpha$ is torsion and $\beta$ is a vector bundle, then $\hom(\alpha,\beta)=0$,
and so we have $\ext(\alpha,\beta)=-\langle\alpha,\beta\rangle$.

Thus suppose that $\alpha$ and $\beta$ are vector bundle
dimension types, with $\alpha$ given by (\ref{e:alphaform})
and $\beta$ by the corresponding formula with coefficients
$m'_*,m'_{ij},d'$.

Suppose that $\eta$ is a dimension type with
$\alpha-\eta\hookrightarrow\alpha$ and $\eta\hookrightarrow \beta$.
Let $\eta$ be given by (\ref{e:alphaform}) with coefficients $m''_*,m''_{ij},d''$.
Clearly
\begin{equation}\label{e:etamin}
m''_* \le \min\{m_*,m'_*\}.
\end{equation}
Since $\eta\hookrightarrow \beta$, it must be a vector bundle
dimension type, so
\begin{equation}\label{e:etapb}
m''_* \ge m''_{i1} \ge m''_{i2} \ge \dots \ge m''_{i,p_i-1} \ge 0
\end{equation}
for all $i$, $m''_*\le m'_*$ and $d''=0$ if $m''_*=0$.

The general sheaf $X$ of dimension $\alpha$ is a vector bundle,
has a sub-bundle $Y$ of dimension type $\alpha-\eta$,
and has underlying vector bundle given by Lemma~\ref{l:underbun}.
Then the underlying bundle of $Y$ is a direct sum of line bundles
of the form $\OO(s)$ with $s\le t+1$, so that
\begin{equation}\label{e:etaeq1}
d-d'' \le m''_*(t+1),
\end{equation}
where $t = \Floor(d/m_*)$.
Dually, the general sheaf of dimension $\beta$ is a vector bundle,
has a sub-bundle of dimension type $\eta$, and has underlying
vector bundle given by Lemma~\ref{l:underbun}. The same argument
now shows that
\begin{equation}\label{e:etaeq2}
d'' \le m''_*(t'+1)
\end{equation}
where $t' = \Floor(d'/m'_*)$.

Clearly there are only finitely many possible $\eta$
satisfying (\ref{e:etamin})-(\ref{e:etaeq2}). For any
such $\eta$, assuming that $\eta,\alpha-\eta\in K_0(\CC)_+$,
one can determine whether or not $\alpha-\eta\hookrightarrow\alpha$
and $\eta\hookrightarrow \beta$ by computing
$\ext(\alpha-\eta,\eta)=0$ and $\ext(\eta,\beta-\eta)=0$.
These are known by the recursion on the rank of $\alpha+\beta$.
\end{proof}

\end{document}